\documentclass{gtart}

\def\ifplaintex{\expandafter\ifx\csname documentclass\endcsname\relax}

\def\gtp{{\mathsurround=0pt\it $\cal G\mskip-2mu$eometry \&\ 
$\cal T\!\!$opology $\cal P\!$ublications}}  

\def\recd{{\small Received:\qua\receiveddate\ifx\reviseddate\relax
\else\qquad Revised:\qua\reviseddate\fi\par}} 

\def\lognumber#1{\def\thelognumber{#1}}
\def\volumenumber#1{\def\thevolumenumber{#1}}
\def\volumeyear#1{\def\thevolumeyear{#1}}
\def\papernumber#1{\def\thepapernumber{#1}}
\def\pagenumbers#1#2{\def\startpage{#1}\def\finishpage{#2}}
\def\published#1{\def\publishdate{#1}}

\def\received#1{\def\receiveddate{#1}}

\def\accepted#1{\def\accepteddate{#1}}

\let\\\par\let\thelognumber\relax\let\thevolumenumber\relax
\let\thepapernumber\relax\let\thevolumeyear\relax\let\startpage\relax
\let\finishpage\relax\let\publishdate\relax\let\receiveddate\relax
\let\reviseddate\relax\let\accepteddate\relax\let\theasciititle\relax
\let\theasciiauthors\relax
\let\theasciiabstract\relax

\let\theasciiemail\relax

\ifplaintex
\font\logobig=cmssbx10 scaled 3836
\font\logomed=cmssbx10 scaled 2557
\else
\font\logobig=cmssbx10 scaled 4200
\font\logomed=cmssbx10 scaled 2800
\fi

\long\def\makeagttitle{   
\count0=\startpage
\agt\hfill      
\hbox to 45truept{\vbox to 0pt{\vglue -13truept{\logomed A\kern -.37em{\logobig 
T}\kern -.38em G}\vss}\hss}
\break
{\small Volume \thevolumenumber\ (\thevolumeyear)
\startpage--\finishpage\nl
Published: \publishdate}

\vglue .25truein

{\parskip=0pt\leftskip 0pt plus
1fil\def\\{\par\smallskip}{\Large\bf\thetitle}\par\medskip} \vglue
0.05truein

{\parskip=0pt\leftskip 0pt plus 1fil\def\\{\par}{\sc\theauthors}
\par\medskip}%
 
\vglue 0.03truein 

{\small\leftskip 25truept\rightskip 25truept{\bf Abstract}\stdspace\theabstract

{\bf AMS Classification}\stdspace\theprimaryclass
\ifx\thesecondaryclass\relax\else; \thesecondaryclass\fi\par
{\bf Keywords}\stdspace \thekeywords\par}\vglue 7truept

}   

\ifplaintex
\font\phead=cmsl9 scaled 950
\font\pnum=cmbx10 scaled 913
\font\pfoot=cmsl9 scaled 950
\headline{\vbox to 0pt{\vskip -4.5mm\line{\small\phead\ifnum
\count0=\startpage ISSN 1472-2739 (on-line) 1472-2747 (printed)
\hfill {\pnum\folio}\else\ifodd\count0\def\\{ }%
\ifx\theshorttitle\relax\thetitle\else\theshorttitle\fi\hfill{\pnum\folio}
\else\def\\{ and }{\pnum\folio}\hfill\ifx\theshortauthors\relax\theauthors
\else\theshortauthors\fi\fi\fi}\vss}}
\footline{\vbox to 0pt{\vglue 0mm\line{\small\pfoot\ifnum\count0=\startpage
\copyright\ \gtp\hfill\else
\agt, Volume \thevolumenumber\ (\thevolumeyear)\hfill\fi}\vss}}
\else
\font=cmsl9 scaled 1050
\font\lnum=cmbx10 
\font=cmsl9 scaled 1050
\makeatletter
\def\@oddhead{{\small\lhead\ifnum\count0=\startpage ISSN 1472-2739 
(on-line) 1472-2747 (printed)\hfill {\lnum\number\count0}\else\ifodd\count0
\def\\{ }\ifx\theshorttitle\relax \thetitle \else\theshorttitle\fi\hfill
{\lnum\number\count0}\else\def\\{ and }{\lnum\number\count0}
\hfill\ifx\theshortauthors\relax 
\theauthors\else\theshortauthors\fi\fi\fi}}\def\@evenhead{\@oddhead}
\def\@oddfoot{\small\lfoot\ifnum\count0=\startpage\copyright\ \gtp\hfill\else
\agt, Volume \thevolumenumber\ (\thevolumeyear)\hfill\fi}
\def\@evenfoot{\@oddfoot}
\makeatother
\fi
\let\maketitlepage\makeagttitle
\let\makeshorttitle\maketitlepage
\let\maketitle\maketitlepage

\newwrite\gtoutfile
\long\gdef\makeheadfile{  
{\def\\{, }\def\s{ }
\immediate\openout\gtoutfile head.xxx
\immediate\write\gtoutfile{To: math@arxiv.org}
\immediate\write\gtoutfile{Subject: put OR rep NNNNN:ppppp}
\immediate\write\gtoutfile{--text follows this line--}
\immediate\write\gtoutfile{Proxy-for: \ifx\theasciiauthors\relax
\theauthors\else\theasciiauthors\fi\s<\ifx\theasciiemail\relax\theemail\else\theasciiemail\fi>}
\immediate\write\gtoutfile{\noexpand\\}
\immediate\write\gtoutfile{Authors: \ifx\theasciiauthors\relax
\theauthors\else\theasciiauthors\fi}
{\def\\{ }\immediate\write\gtoutfile{Title: \ifx\theasciititle\relax
\thetitle\else\theasciititle\fi}}
\immediate\write\gtoutfile{Subj-class: GT or SG, GR etc}
\immediate\write\gtoutfile{MSC-class: \theprimaryclass\ifx\thesecondaryclass\relax\else, \thesecondaryclass\fi}
\immediate\write\gtoutfile{Journal-ref: Algebr. Geom. Topol. \thevolumenumber\s
(\thevolumeyear) \startpage-\finishpage}
\immediate\write\gtoutfile{Comments: Published by Algebraic and
Geometric Topology at}
\immediate\write\gtoutfile{\s\s\s  http://www.maths.warwick.ac.uk/agt/AGTVol\thevolumenumber/agt-\thevolumenumber-\thepapernumber.abs.html}
\immediate\write\gtoutfile{\noexpand\\}
\immediate\write\gtoutfile{}
\ifx\theasciiabstract\relax
\immediate\write\gtoutfile{\theabstract}\else
\immediate\write\gtoutfile{\theasciiabstract}\fi
\immediate\write\gtoutfile{}
\immediate\write\gtoutfile{\noexpand\\}
\immediate\write\gtoutfile{}
\immediate\closeout\gtoutfile}}  

\def\maketitlepage{\makeagttitle\makeheadfile}
\let\makeshorttitle\maketitlepage
\let\maketitle\maketitlepage

\def\ifplaintex{\expandafter\ifx\csname documentclass\endcsname\relax}

\def\gtp{{\mathsurround=0pt\it $\cal G\mskip-2mu$eometry \&\ 
$\cal T\!\!$opology $\cal P\!$ublications}}  

\def\recd{{\small Received:\qua\receiveddate\ifx\reviseddate\relax
\else\qquad Revised:\qua\reviseddate\fi\par}} 

\def\lognumber#1{\def\thelognumber{#1}}
\def\volumenumber#1{\def\thevolumenumber{#1}}
\def\volumeyear#1{\def\thevolumeyear{#1}}
\def\papernumber#1{\def\thepapernumber{#1}}
\def\pagenumbers#1#2{\def\startpage{#1}\def\finishpage{#2}}
\def\published#1{\def\publishdate{#1}}

\def\received#1{\def\receiveddate{#1}}

\def\accepted#1{\def\accepteddate{#1}}

\let\\\par\let\thelognumber\relax\let\thevolumenumber\relax
\let\thepapernumber\relax\let\thevolumeyear\relax\let\startpage\relax
\let\finishpage\relax\let\publishdate\relax\let\receiveddate\relax
\let\reviseddate\relax\let\accepteddate\relax\let\theasciititle\relax
\let\theasciiauthors\relax
\let\theasciiabstract\relax

\let\theasciiemail\relax

\ifplaintex
\font\logobig=cmssbx10 scaled 3836
\font\logomed=cmssbx10 scaled 2557
\else
\font\logobig=cmssbx10 scaled 4200
\font\logomed=cmssbx10 scaled 2800
\fi

\long\def\makeagttitle{   
\count0=\startpage
\agt\hfill      
\hbox to 45truept{\vbox to 0pt{\vglue -13truept{\logomed A\kern -.37em{\logobig 
T}\kern -.38em G}\vss}\hss}
\break
{\small Volume \thevolumenumber\ (\thevolumeyear)
\startpage--\finishpage\nl
Published: \publishdate}

\vglue .25truein

{\parskip=0pt\leftskip 0pt plus
1fil\def\\{\par\smallskip}{\Large\bf\thetitle}\par\medskip} \vglue
0.05truein

{\parskip=0pt\leftskip 0pt plus 1fil\def\\{\par}{\sc\theauthors}
\par\medskip}%
 
\vglue 0.03truein 

{\small\leftskip 25truept\rightskip 25truept{\bf Abstract}\stdspace\theabstract

{\bf AMS Classification}\stdspace\theprimaryclass
\ifx\thesecondaryclass\relax\else; \thesecondaryclass\fi\par
{\bf Keywords}\stdspace \thekeywords\par}\vglue 7truept

}   

\ifplaintex
\font\phead=cmsl9 scaled 950
\font\pnum=cmbx10 scaled 913
\font\pfoot=cmsl9 scaled 950
\headline{\vbox to 0pt{\vskip -4.5mm\line{\small\phead\ifnum
\count0=\startpage ISSN 1472-2739 (on-line) 1472-2747 (printed)
\hfill {\pnum\folio}\else\ifodd\count0\def\\{ }%
\ifx\theshorttitle\relax\thetitle\else\theshorttitle\fi\hfill{\pnum\folio}
\else\def\\{ and }{\pnum\folio}\hfill\ifx\theshortauthors\relax\theauthors
\else\theshortauthors\fi\fi\fi}\vss}}
\footline{\vbox to 0pt{\vglue 0mm\line{\small\pfoot\ifnum\count0=\startpage
\copyright\ \gtp\hfill\else
\agt, Volume \thevolumenumber\ (\thevolumeyear)\hfill\fi}\vss}}
\else
\font=cmsl9 scaled 1050
\font\lnum=cmbx10 
\font=cmsl9 scaled 1050
\makeatletter
\def\@oddhead{{\small\lhead\ifnum\count0=\startpage ISSN 1472-2739 
(on-line) 1472-2747 (printed)\hfill {\lnum\number\count0}\else\ifodd\count0
\def\\{ }\ifx\theshorttitle\relax \thetitle \else\theshorttitle\fi\hfill
{\lnum\number\count0}\else\def\\{ and }{\lnum\number\count0}
\hfill\ifx\theshortauthors\relax 
\theauthors\else\theshortauthors\fi\fi\fi}}\def\@evenhead{\@oddhead}
\def\@oddfoot{\small\lfoot\ifnum\count0=\startpage\copyright\ \gtp\hfill\else
\agt, Volume \thevolumenumber\ (\thevolumeyear)\hfill\fi}
\def\@evenfoot{\@oddfoot}
\makeatother
\fi
\let\maketitlepage\makeagttitle
\let\makeshorttitle\maketitlepage
\let\maketitle\maketitlepage

\newwrite\gtoutfile
\long\gdef\makeheadfile{  
{\def\\{, }\def\s{ }
\immediate\openout\gtoutfile head.xxx
\immediate\write\gtoutfile{To: math@arxiv.org}
\immediate\write\gtoutfile{Subject: put OR rep NNNNN:ppppp}
\immediate\write\gtoutfile{--text follows this line--}
\immediate\write\gtoutfile{Proxy-for: \ifx\theasciiauthors\relax
\theauthors\else\theasciiauthors\fi\s<\ifx\theasciiemail\relax\theemail\else\theasciiemail\fi>}
\immediate\write\gtoutfile{\noexpand\\}
\immediate\write\gtoutfile{Authors: \ifx\theasciiauthors\relax
\theauthors\else\theasciiauthors\fi}
{\def\\{ }\immediate\write\gtoutfile{Title: \ifx\theasciititle\relax
\thetitle\else\theasciititle\fi}}
\immediate\write\gtoutfile{Subj-class: GT or SG, GR etc}
\immediate\write\gtoutfile{MSC-class: \theprimaryclass\ifx\thesecondaryclass\relax\else, \thesecondaryclass\fi}
\immediate\write\gtoutfile{Journal-ref: Algebr. Geom. Topol. \thevolumenumber\s
(\thevolumeyear) \startpage-\finishpage}
\immediate\write\gtoutfile{Comments: Published by Algebraic and
Geometric Topology at}
\immediate\write\gtoutfile{\s\s\s  http://www.maths.warwick.ac.uk/agt/AGTVol\thevolumenumber/agt-\thevolumenumber-\thepapernumber.abs.html}
\immediate\write\gtoutfile{\noexpand\\}
\immediate\write\gtoutfile{}
\ifx\theasciiabstract\relax
\immediate\write\gtoutfile{\theabstract}\else
\immediate\write\gtoutfile{\theasciiabstract}\fi
\immediate\write\gtoutfile{}
\immediate\write\gtoutfile{\noexpand\\}
\immediate\write\gtoutfile{}
\immediate\closeout\gtoutfile}}  

\def\maketitlepage{\makeagttitle\makeheadfile}
\let\makeshorttitle\maketitlepage
\let\maketitle\maketitlepage

\def\ifplaintex{\expandafter\ifx\csname documentclass\endcsname\relax}

\def\gtp{{\mathsurround=0pt\it $\cal G\mskip-2mu$eometry \&\ 
$\cal T\!\!$opology $\cal P\!$ublications}}  

\def\recd{{\small Received:\qua\receiveddate\ifx\reviseddate\relax
\else\qquad Revised:\qua\reviseddate\fi\par}} 

\def\lognumber#1{\def\thelognumber{#1}}
\def\volumenumber#1{\def\thevolumenumber{#1}}
\def\volumeyear#1{\def\thevolumeyear{#1}}
\def\papernumber#1{\def\thepapernumber{#1}}
\def\pagenumbers#1#2{\def\startpage{#1}\def\finishpage{#2}}
\def\published#1{\def\publishdate{#1}}

\def\received#1{\def\receiveddate{#1}}

\def\accepted#1{\def\accepteddate{#1}}

\let\\\par\let\thelognumber\relax\let\thevolumenumber\relax
\let\thepapernumber\relax\let\thevolumeyear\relax\let\startpage\relax
\let\finishpage\relax\let\publishdate\relax\let\receiveddate\relax
\let\reviseddate\relax\let\accepteddate\relax\let\theasciititle\relax
\let\theasciiauthors\relax
\let\theasciiabstract\relax

\let\theasciiemail\relax

\ifplaintex
\font\logobig=cmssbx10 scaled 3836
\font\logomed=cmssbx10 scaled 2557
\else
\font\logobig=cmssbx10 scaled 4200
\font\logomed=cmssbx10 scaled 2800
\fi

\long\def\makeagttitle{   
\count0=\startpage
\agt\hfill      
\hbox to 45truept{\vbox to 0pt{\vglue -13truept{\logomed A\kern -.37em{\logobig 
T}\kern -.38em G}\vss}\hss}
\break
{\small Volume \thevolumenumber\ (\thevolumeyear)
\startpage--\finishpage\nl
Published: \publishdate}

\vglue .25truein

{\parskip=0pt\leftskip 0pt plus
1fil\def\\{\par\smallskip}{\Large\bf\thetitle}\par\medskip} \vglue
0.05truein

{\parskip=0pt\leftskip 0pt plus 1fil\def\\{\par}{\sc\theauthors}
\par\medskip}%
 
\vglue 0.03truein 

{\small\leftskip 25truept\rightskip 25truept{\bf Abstract}\stdspace\theabstract

{\bf AMS Classification}\stdspace\theprimaryclass
\ifx\thesecondaryclass\relax\else; \thesecondaryclass\fi\par
{\bf Keywords}\stdspace \thekeywords\par}\vglue 7truept

}   

\ifplaintex
\font\phead=cmsl9 scaled 950
\font\pnum=cmbx10 scaled 913
\font\pfoot=cmsl9 scaled 950
\headline{\vbox to 0pt{\vskip -4.5mm\line{\small\phead\ifnum
\count0=\startpage ISSN 1472-2739 (on-line) 1472-2747 (printed)
\hfill {\pnum\folio}\else\ifodd\count0\def\\{ }%
\ifx\theshorttitle\relax\thetitle\else\theshorttitle\fi\hfill{\pnum\folio}
\else\def\\{ and }{\pnum\folio}\hfill\ifx\theshortauthors\relax\theauthors
\else\theshortauthors\fi\fi\fi}\vss}}
\footline{\vbox to 0pt{\vglue 0mm\line{\small\pfoot\ifnum\count0=\startpage
\copyright\ \gtp\hfill\else
\agt, Volume \thevolumenumber\ (\thevolumeyear)\hfill\fi}\vss}}
\else
\font=cmsl9 scaled 1050
\font\lnum=cmbx10 
\font=cmsl9 scaled 1050
\makeatletter
\def\@oddhead{{\small\lhead\ifnum\count0=\startpage ISSN 1472-2739 
(on-line) 1472-2747 (printed)\hfill {\lnum\number\count0}\else\ifodd\count0
\def\\{ }\ifx\theshorttitle\relax \thetitle \else\theshorttitle\fi\hfill
{\lnum\number\count0}\else\def\\{ and }{\lnum\number\count0}
\hfill\ifx\theshortauthors\relax 
\theauthors\else\theshortauthors\fi\fi\fi}}\def\@evenhead{\@oddhead}
\def\@oddfoot{\small\lfoot\ifnum\count0=\startpage\copyright\ \gtp\hfill\else
\agt, Volume \thevolumenumber\ (\thevolumeyear)\hfill\fi}
\def\@evenfoot{\@oddfoot}
\makeatother
\fi
\let\maketitlepage\makeagttitle
\let\makeshorttitle\maketitlepage
\let\maketitle\maketitlepage

\newwrite\gtoutfile
\long\gdef\makeheadfile{  
{\def\\{, }\def\s{ }
\immediate\openout\gtoutfile head.xxx
\immediate\write\gtoutfile{To: math@arxiv.org}
\immediate\write\gtoutfile{Subject: put OR rep NNNNN:ppppp}
\immediate\write\gtoutfile{--text follows this line--}
\immediate\write\gtoutfile{Proxy-for: \ifx\theasciiauthors\relax
\theauthors\else\theasciiauthors\fi\s<\ifx\theasciiemail\relax\theemail\else\theasciiemail\fi>}
\immediate\write\gtoutfile{\noexpand\\}
\immediate\write\gtoutfile{Authors: \ifx\theasciiauthors\relax
\theauthors\else\theasciiauthors\fi}
{\def\\{ }\immediate\write\gtoutfile{Title: \ifx\theasciititle\relax
\thetitle\else\theasciititle\fi}}
\immediate\write\gtoutfile{Subj-class: GT or SG, GR etc}
\immediate\write\gtoutfile{MSC-class: \theprimaryclass\ifx\thesecondaryclass\relax\else, \thesecondaryclass\fi}
\immediate\write\gtoutfile{Journal-ref: Algebr. Geom. Topol. \thevolumenumber\s
(\thevolumeyear) \startpage-\finishpage}
\immediate\write\gtoutfile{Comments: Published by Algebraic and
Geometric Topology at}
\immediate\write\gtoutfile{\s\s\s  http://www.maths.warwick.ac.uk/agt/AGTVol\thevolumenumber/agt-\thevolumenumber-\thepapernumber.abs.html}
\immediate\write\gtoutfile{\noexpand\\}
\immediate\write\gtoutfile{}
\ifx\theasciiabstract\relax
\immediate\write\gtoutfile{\theabstract}\else
\immediate\write\gtoutfile{\theasciiabstract}\fi
\immediate\write\gtoutfile{}
\immediate\write\gtoutfile{\noexpand\\}
\immediate\write\gtoutfile{}
\immediate\closeout\gtoutfile}}  

\def\maketitlepage{\makeagttitle\makeheadfile}
\let\makeshorttitle\maketitlepage
\let\maketitle\maketitlepage

\def\ifplaintex{\expandafter\ifx\csname documentclass\endcsname\relax}

\def\gtp{{\mathsurround=0pt\it $\cal G\mskip-2mu$eometry \&\ 
$\cal T\!\!$opology $\cal P\!$ublications}}  

\def\recd{{\small Received:\qua\receiveddate\ifx\reviseddate\relax
\else\qquad Revised:\qua\reviseddate\fi\par}} 

\def\lognumber#1{\def\thelognumber{#1}}
\def\volumenumber#1{\def\thevolumenumber{#1}}
\def\volumeyear#1{\def\thevolumeyear{#1}}
\def\papernumber#1{\def\thepapernumber{#1}}
\def\pagenumbers#1#2{\def\startpage{#1}\def\finishpage{#2}}
\def\published#1{\def\publishdate{#1}}

\def\received#1{\def\receiveddate{#1}}

\def\accepted#1{\def\accepteddate{#1}}

\let\\\par\let\thelognumber\relax\let\thevolumenumber\relax
\let\thepapernumber\relax\let\thevolumeyear\relax\let\startpage\relax
\let\finishpage\relax\let\publishdate\relax\let\receiveddate\relax
\let\reviseddate\relax\let\accepteddate\relax\let\theasciititle\relax
\let\theasciiauthors\relax
\let\theasciiabstract\relax

\let\theasciiemail\relax

\ifplaintex
\font\logobig=cmssbx10 scaled 3836
\font\logomed=cmssbx10 scaled 2557
\else
\font\logobig=cmssbx10 scaled 4200
\font\logomed=cmssbx10 scaled 2800
\fi

\long\def\makeagttitle{   
\count0=\startpage
\agt\hfill      
\hbox to 45truept{\vbox to 0pt{\vglue -13truept{\logomed A\kern -.37em{\logobig 
T}\kern -.38em G}\vss}\hss}
\break
{\small Volume \thevolumenumber\ (\thevolumeyear)
\startpage--\finishpage\nl
Published: \publishdate}

\vglue .25truein

{\parskip=0pt\leftskip 0pt plus
1fil\def\\{\par\smallskip}{\Large\bf\thetitle}\par\medskip} \vglue
0.05truein

{\parskip=0pt\leftskip 0pt plus 1fil\def\\{\par}{\sc\theauthors}
\par\medskip}%
 
\vglue 0.03truein 

{\small\leftskip 25truept\rightskip 25truept{\bf Abstract}\stdspace\theabstract

{\bf AMS Classification}\stdspace\theprimaryclass
\ifx\thesecondaryclass\relax\else; \thesecondaryclass\fi\par
{\bf Keywords}\stdspace \thekeywords\par}\vglue 7truept

}   

\ifplaintex
\font\phead=cmsl9 scaled 950
\font\pnum=cmbx10 scaled 913
\font\pfoot=cmsl9 scaled 950
\headline{\vbox to 0pt{\vskip -4.5mm\line{\small\phead\ifnum
\count0=\startpage ISSN 1472-2739 (on-line) 1472-2747 (printed)
\hfill {\pnum\folio}\else\ifodd\count0\def\\{ }%
\ifx\theshorttitle\relax\thetitle\else\theshorttitle\fi\hfill{\pnum\folio}
\else\def\\{ and }{\pnum\folio}\hfill\ifx\theshortauthors\relax\theauthors
\else\theshortauthors\fi\fi\fi}\vss}}
\footline{\vbox to 0pt{\vglue 0mm\line{\small\pfoot\ifnum\count0=\startpage
\copyright\ \gtp\hfill\else
\agt, Volume \thevolumenumber\ (\thevolumeyear)\hfill\fi}\vss}}
\else
\font=cmsl9 scaled 1050
\font\lnum=cmbx10 
\font=cmsl9 scaled 1050
\makeatletter
\def\@oddhead{{\small\lhead\ifnum\count0=\startpage ISSN 1472-2739 
(on-line) 1472-2747 (printed)\hfill {\lnum\number\count0}\else\ifodd\count0
\def\\{ }\ifx\theshorttitle\relax \thetitle \else\theshorttitle\fi\hfill
{\lnum\number\count0}\else\def\\{ and }{\lnum\number\count0}
\hfill\ifx\theshortauthors\relax 
\theauthors\else\theshortauthors\fi\fi\fi}}\def\@evenhead{\@oddhead}
\def\@oddfoot{\small\lfoot\ifnum\count0=\startpage\copyright\ \gtp\hfill\else
\agt, Volume \thevolumenumber\ (\thevolumeyear)\hfill\fi}
\def\@evenfoot{\@oddfoot}
\makeatother
\fi
\let\maketitlepage\makeagttitle
\let\makeshorttitle\maketitlepage
\let\maketitle\maketitlepage

\newwrite\gtoutfile
\long\gdef\makeheadfile{  
{\def\\{, }\def\s{ }
\immediate\openout\gtoutfile head.xxx
\immediate\write\gtoutfile{To: math@arxiv.org}
\immediate\write\gtoutfile{Subject: put OR rep NNNNN:ppppp}
\immediate\write\gtoutfile{--text follows this line--}
\immediate\write\gtoutfile{Proxy-for: \ifx\theasciiauthors\relax
\theauthors\else\theasciiauthors\fi\s<\ifx\theasciiemail\relax\theemail\else\theasciiemail\fi>}
\immediate\write\gtoutfile{\noexpand\\}
\immediate\write\gtoutfile{Authors: \ifx\theasciiauthors\relax
\theauthors\else\theasciiauthors\fi}
{\def\\{ }\immediate\write\gtoutfile{Title: \ifx\theasciititle\relax
\thetitle\else\theasciititle\fi}}
\immediate\write\gtoutfile{Subj-class: GT or SG, GR etc}
\immediate\write\gtoutfile{MSC-class: \theprimaryclass\ifx\thesecondaryclass\relax\else, \thesecondaryclass\fi}
\immediate\write\gtoutfile{Journal-ref: Algebr. Geom. Topol. \thevolumenumber\s
(\thevolumeyear) \startpage-\finishpage}
\immediate\write\gtoutfile{Comments: Published by Algebraic and
Geometric Topology at}
\immediate\write\gtoutfile{\s\s\s  http://www.maths.warwick.ac.uk/agt/AGTVol\thevolumenumber/agt-\thevolumenumber-\thepapernumber.abs.html}
\immediate\write\gtoutfile{\noexpand\\}
\immediate\write\gtoutfile{}
\ifx\theasciiabstract\relax
\immediate\write\gtoutfile{\theabstract}\else
\immediate\write\gtoutfile{\theasciiabstract}\fi
\immediate\write\gtoutfile{}
\immediate\write\gtoutfile{\noexpand\\}
\immediate\write\gtoutfile{}
\immediate\closeout\gtoutfile}}  

\def\maketitlepage{\makeagttitle\makeheadfile}
\let\makeshorttitle\maketitlepage
\let\maketitle\maketitlepage

\lognumber{6}
\volumenumber{2}
\volumeyear{2002}
\papernumber{6}
\published{14 February 2002}
\pagenumbers{95}{135}
\received{27 September 2001}
\accepted{8 February 2002}

\usepackage{amssymb,epsf,amsmath,amscd}

\newtheorem{theorem}{Theorem}[section]
\newtheorem{corollary}[theorem]{Corollary}
\newtheorem{lemma}[theorem]{Lemma}
\newtheorem{proposition}[theorem]{Proposition}

\theoremstyle{definition}
\newtheorem{definition}[theorem]{Definition}
\newtheorem{example}[theorem]{Example}
\newtheorem{remark}[theorem]{Remark}

\newcommand{\Z}{\mathbb{Z}}
\newcommand{\R}{\mathbb{R}}
\newcommand{\LX}{ {\cal L} }

\begin{document}

\title{Twisted quandle homology theory\\and cocycle knot invariants}

\authors{J. Scott Carter\\Mohamed Elhamdadi\\Masahico Saito}
\shortauthors{Carter, Elhamdadi and Saito}

\address{University of South Alabama, Mobile, AL 36688, USA\\
University of South Florida, Tampa, FL 33620, USA\\
University of South Florida, Tampa, FL 33620, USA}  

\emails{carter@mathstat.usouthal.edu, emohamed@math.usf.edu,
saito@math.usf.edu}

\begin{abstract}
The quandle homology theory is generalized to the case when 
the coefficient groups admit  the structure of Alexander quandles, 
by including  an action of the infinite cyclic group in the boundary operator.
Theories of Alexander extensions of quandles in relation
to low dimensional cocycles are developed in parallel to
group extension theories for group cocycles. Explicit formulas for cocycles 
corresponding to extensions are given, and used to prove
non-triviality of cohomology groups for some quandles. 
The corresponding generalization of the quandle cocycle knot invariants
is given, by using the Alexander numbering of regions 
in the definition of state-sums. The invariants are used to derive 
information on twisted cohomology groups.
\end{abstract}


\maketitle

\section{Introduction}

A quandle is a set with a 
self-distributive binary operation (defined below)
whose definition was partially motivated from knot theory. 
A (co)homology theory was defined in \cite{CJKLS} for quandles,
which is a modification of rack (co)homology defined in \cite{FRS2}. 
State-sum invariants,
called the quandle cocycle invariants, 
 using quandle cocycles as weights are 
defined \cite{CJKLS} and computed for important families
of classical knots and knotted surfaces \cite{CJKS1}.
Quandle homomorphisms and virtual knots are applied to this 
homology theory \cite{CJKS2}.
 The invariants were applied to study 
knots, 
for example, in detecting non-invertible
knotted surfaces \cite{CJKLS}. 
On the other hand, knot diagrams colored by quandles can be used 
to study quandle homology groups. This viewpoint was developed
 in \cite{FRS2,Flower,Greene}  
for rack homology and homotopy, and generalized to quandle homology
in \cite{SSS2}.
Thus, the algebraic theory of quandle homology has been applied to 
knot invariants, and geometric methods using knot diagrams
have been applied to quandle homology theory.

Computations of (co)homology groups, however, had  depended upon
computer assisted calculations, until in \cite{SanFran},
relations of low dimensional cocycles to extensions of quandles 
were given.
These were 
used in \cite{ext} 
to give an algebraic method of constructing cocycles
explicitly and to obtain new cocycles via quandle extensions.
The methods introduced in \cite{ext} are developed to
parallel the theory of group $2$-cocycles in relation to group
extensions \cite{Brown}. 

In this paper, we develop the method of quandle extensions
when the coefficient group admits the structure of a 
${\Z}[T,T^{-1}]$-module.
 In this case, the coefficients also have a quandle structure
 and new cocycles arise via the theory of extensions. 
This theory of twisted coefficients is an analogue of group and 
Hoshschild cohomology in which the coefficient rings admit
actions. State-sum invariants can be obtained from the twisted
 cohomology theory using Alexander numbering on the regions of
 the knot diagram. These invariants then yield information on the twisted 
quandle cohomology groups.

The paper is organized as follows. 
In Section~\ref{revsec}, necessary materials are reviewed briefly.
The twisted quandle homology theory is defined in Section~\ref{defsec}, 
and a few examples are given.
The obstruction and extension theories are developed for low dimensional
cocycles in Section~\ref{extsec}, and families of Alexander quandles
are presented in Section~\ref{lxsec} as examples. 
Explicit  
formulas for cocycles are also provided.
In Section~\ref{H1coeffsec}, cohomology groups with cohomology coefficients
are used for further constructions of cocycles. 
In Section~\ref{invsec}, the twisted cocycles are used to generalize
cocycle knot invariants, using Alexander numbering of regions,
and applications are given.

\rk{Acknowledgements} JSC was supported in part by NSF Grant 
DMS\,9988107.  MS was supported in part by NSF Grant DMS 9988101.  The
authors would like to thank the referee for carefully reading the
manuscript and suggesting improvements.

\section{Quandles and their homology theory} \label{revsec}

In this section we review necessary material from the papers mentioned
in the introduction.

A {\it quandle}, $X$, is a set with a binary operation $(a, b) \mapsto a * b$
such that

(I) For any $a \in X$,
$a* a =a$.

(II) For any $a,b \in X$, there is a unique $c \in X$ such that 
$a= c*b$.

(III) 
For any $a,b,c \in X$, we have
$ (a*b)*c=(a*c)*(b*c). $

A {\it rack} is a set with a binary operation that satisfies 
(II) and (III).
Racks and quandles have been studied in, for example, 
\cite{Brieskorn,FR,Joyce,K&P,Matveev}.
The axioms for a quandle correspond respectively to the 
Reidemeister moves of type I, II, and III 
(see 
\cite{FR,K&P}, for example). 
A function $f: X \rightarrow  Y$ between quandles
or racks  is a {\it homomorphism}
if $f(a \ast b) = f(a) * f(b)$ 
for any $a, b \in X$.

The following are typical examples of quandles.

\begin{itemize}
\item
A group $X=G$ with
$n$-fold conjugation
as the quandle operation: $a*b=b^{-n} a b^n$.
\item
Any set $X$ with the operation $x*y=x$ for any $x,y \in X$ is
a quandle called the {\it trivial} quandle.
The trivial quandle of $n$ elements is denoted by $T_n$.
\item
Let $n$ be a positive integer.
For elements  $i, j \in \{ 0, 1, \ldots , n-1 \}$, define
$i\ast j \equiv 2j-i \pmod{n}$.
Then $\ast$ defines a quandle
structure  called the {\it dihedral quandle},
  $R_n$.
This set can be identified with  the
set of reflections of a regular $n$-gon
  with conjugation
as the quandle operation.
\item
Any $\Lambda (={\Z}[T, T^{-1}])$-module $M$
is a quandle with
$a*b=Ta+(1-T)b$, $a,b \in M$, called an {\it  Alexander  quandle}.
Furthermore for a positive integer
$n$, a {\it mod-$n$ Alexander  quandle}
${\Z}_n[T, T^{-1}]/(h(T))$
is a quandle
for
a Laurent polynomial $h(T)$.
The mod-$n$ Alexander quandle is finite
if the coefficients of the
highest and lowest degree terms
of $h$
  are units of ${\Z}_n$.
\end{itemize}

\begin{figure}[ht!]
\begin{center}
\mbox{
\epsfxsize=4in
\epsfbox{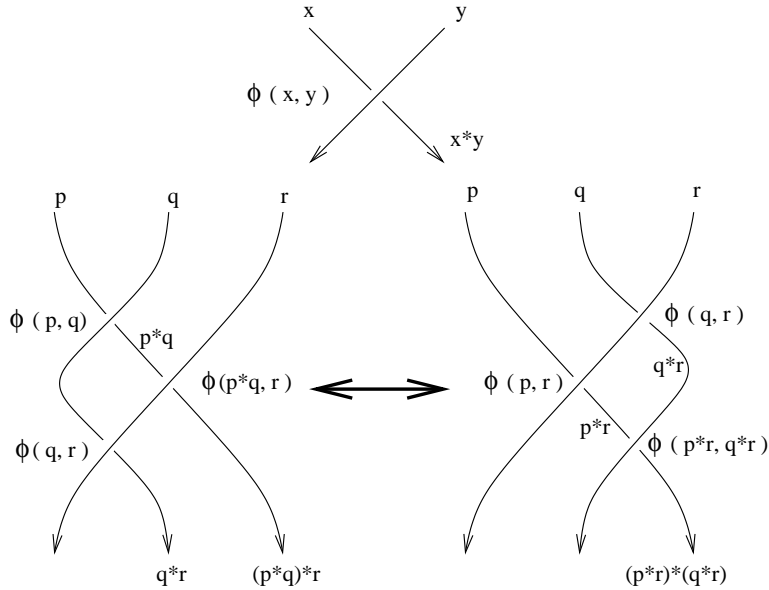} 
}
\end{center}
\caption{ Type III move and the quandle identity  }
\label{2cocy} 
\end{figure}

Let $C_n^{\rm R}(X)$ be the free 
abelian group generated by
$n$-tuples $(x_1, \dots, x_n)$ of elements of a quandle $X$. Define a
homomorphism
$\partial_{n}: C_{n}^{\rm R}(X) \to C_{n-1}^{\rm R}(X)$ by \begin{eqnarray}
\lefteqn{
\partial_{n}(x_1, x_2, \dots, x_n) } \nonumber \\ && =
\sum_{i=2}^{n} (-1)^{i}\left[ (x_1, x_2, \dots, x_{i-1}, x_{i+1},\dots, x_n) \right.
\nonumber \\
&&
- \left. (x_1 \ast x_i, x_2 \ast x_i, \dots, x_{i-1}\ast x_i, x_{i+1}, \dots, x_n) \right]
\end{eqnarray}
for $n \geq 2$ 
and $\partial_n=0$ for 
$n \leq 1$. 
 Then
$C_\ast^{\rm R}(X)
= \{C_n^{\rm R}(X), \partial_n \}$ is a chain complex.

Let $C_n^{\rm D}(X)$ be the subset of $C_n^{\rm R}(X)$ generated
by $n$-tuples $(x_1, \dots, x_n)$
with $x_{i}=x_{i+1}$ for some $i \in \{1, \dots,n-1\}$ if $n \geq 2$;
otherwise let $C_n^{\rm D}(X)=0$. If $X$ is a quandle, then
$\partial_n(C_n^{\rm D}(X)) \subset C_{n-1}^{\rm D}(X)$ and
$C_\ast^{\rm D}(X) = \{ C_n^{\rm D}(X), \partial_n \}$ is a sub-complex of
$C_\ast^{\rm
R}(X)$. Put $C_n^{\rm Q}(X) = C_n^{\rm R}(X)/ C_n^{\rm D}(X)$ and 
$C_\ast^{\rm Q}(X) = \{ C_n^{\rm Q}(X), \partial'_n \}$,
where $\partial'_n$ is the induced homomorphism.
Henceforth, all boundary maps will be denoted by $\partial_n$.

For an abelian group $G$, define the chain and cochain complexes
\begin{eqnarray}
C_\ast^{\rm W}(X;G) = C_\ast^{\rm W}(X) \otimes G, \quad && \partial =
\partial \otimes {\rm id}; \\ C^\ast_{\rm W}(X;G) = {\rm Hom}(C_\ast^{\rm
W}(X), G), \quad
&& \delta= {\rm Hom}(\partial, {\rm id})
\end{eqnarray}
in the usual way, where ${\rm W}$ 
 $={\rm D}$, ${\rm R}$, ${\rm Q}$.

The groups of cycles and boundaries are denoted respectively by 
${\mbox{\rm ker}}(\partial) =Z_n^{\rm W}(X;G) \subset C_n^{\rm W} (X;G)$ and 
${\mbox{\rm Im}}(\partial)=B_n^{\rm W}(X;G) \subset C_n^{\rm W}(X;G)$ 
while the cocycles and coboundaries are denoted respectively by
${\mbox{\rm ker}}(\delta) =Z^n_{\rm W}(X;G) \subset C^n_{\rm W}(X;G)$ and 
${\mbox{\rm Im}}(\partial)=B^n_{\rm W}(X;G) \subset C^n_{\rm W}(X;G).$
 In particular,  a quandle $2$-cocycle is
an element  $\phi \in Z^2_{\rm Q}(X;G)$,  and 
the equalities 
\begin{eqnarray*}
\phi(x,z)+\phi(x*z,y*z)&=&\phi(x*y,z)+\phi(x,y)  \\
 \mbox{and} \quad  \phi(x,x) & = & 0
\end{eqnarray*}
are satisfied for all $x,y,z\in X$.

The $n$\/th {\it quandle homology group\/}  and the $n$\/th
{\it quandle cohomology group\/ } \cite{CJKLS} of a quandle $X$ with coefficient group $G$ are
\begin{eqnarray}
H_n^{\rm Q}(X; G) 
 &= &H_{n}(C_\ast^{\rm Q}(X;G)) \; = \; Z_n^{\rm Q}(X;G)/B_n^{\rm Q}(X;G), 
 \nonumber \\
H^n_{\rm Q}(X; G) 
 &= & H^{n}(C^\ast_{\rm Q}(X;G))\; = \; Z^n_{\rm Q}(X;G)/B^n_{\rm Q}(X;G). \end{eqnarray}
Let a classical knot diagram be given. 
The co-orientation is a family of normal vectors to the knot diagram 
such that the pair (orientation, co-orientation) 
matches
the given (right-handed, or counterclockwise) orientation of the plane.
At a crossing, 
if the pair of the co-orientation 
 of the 
over-arc and  that of the under-arc
matches the (right-hand) orientation of the plane, then the 
crossing is called {\it positive}; otherwise it is {\it negative}. 
Crossings in Fig.~\ref{2cocy} are positive by convention.

A  {\it coloring}   
of an oriented  classical knot diagram is a
function ${\mathcal C} : R \rightarrow X$, where $X$ is a fixed 
quandle
and $R$ is the set of over-arcs in the diagram,
satisfying the  condition
depicted 
in the top
of Fig.~\ref{2cocy}. 
In the figure, a 
crossing with
over-arc, $r$, has color ${\mathcal C}(r)= y \in X$. 
The under-arcs are called $r_1$ and $r_2$ from top to bottom;
the normal (co-orientation) of the over-arc $r$ points from $r_1$ to $r_2$.
Then it is required that 
${\mathcal C}(r_1)= x$ and ${\mathcal C}(r_2)=x*y$.
Observe that a coloring is a quandle homomorphism (${\mathcal C}(x*y)= {\mathcal C}(x) * {\mathcal C}(y)$) from the fundamental quandle of the knot
 (see \cite{Joyce}) to the quandle $X$.

Note that locally the colors do not depend on the 
orientation of the under-arc.
The quandle element ${\mathcal C}(r)$ assigned to an arc $r$ by a coloring 
 ${\mathcal C}$ is called a {\it color} of the arc. 
This definition of colorings on knot diagrams has been known, see 
\cite{FR,FoxTrip} for example. 
Henceforth, all the quandles that are used to color diagrams will be finite.

In Fig.~\ref{2cocy} bottom, the relation between Redemeister type III move
and a quandle axiom (self-distributivity) is indicated. 
In particular, the colors of the bottom right segments before and after
the move correspond to the self-distributivity. 

Let a quandle $X$, and a 
quandle $2$-cocycle $\phi \in Z^2_{\rm Q}(X;A)$ be given.
A {\it (Boltzmann) weight}, $B(\tau, {\cal C})$ 
(that depends on $\phi$),
at a  crossing $\tau$ is defined as follows.
Let  ${\cal  C}$ 
denote a coloring.
Let $r$ be the over-arc at $\tau$, and $r_1$, $r_2$ be 
under-arcs such that the 
normal
to $r$ points from $r_1$ to $r_2$.
Let $x={\cal C}(r_1)$ and $y={\cal C}(r)$. 
Then define $B(\tau, {\cal C})= \phi(x,y)^{\epsilon (\tau)}$,
where 
$\epsilon (\tau)= 1$ or $-1$, if  the sign of $\tau$ 
is positive or negative, respectively.

The {\it partition function}, or a {\it state-sum}, 
is the expression 
$$
\sum_{{\cal C}}  \prod_{\tau}  B( \tau, {\cal C}).
$$
The product is taken over all crossings of the given diagram,
and the sum is taken over all possible colorings.
The values of the partition function 
are  taken to be in  the group ring ${\Z}[A]$ where $A$ is the coefficient 
group 
written multiplicatively. 
The partition function depends on the choice of $2$-cocycle
$\phi$. 
This is proved \cite{CJKLS} to be a knot invariant, called
the {\it (quandle) cocycle invariant}. 
Figure~\ref{2cocy} shows the invariance of the state-sum under
the Reidemeister type III move.

\section{Twisted  quandle homology } \label{defsec}

In this section we generalize the quandle homology theory to 
those with coefficients in Alexander quandles.

Let $\Lambda=\Z [T, T^{-1}]$, and let 
 $C_n^{\rm TR}(X)= C_n^{\rm TR}(X; \Lambda)$ be the free 
module over $\Lambda$   generated by
$n$-tuples $(x_1, \dots, x_n)$ of elements of a quandle $X$. 
Define a homomorphism 
$\partial = \partial^T_{n}: C_{n}^{\rm TR}(X) \to  C_{n-1}^{\rm TR}(X) $ 
by \begin{eqnarray}
\lefteqn{
\partial^T_{n}(x_1, x_2, \dots, x_n) } \nonumber \\ && =
\sum_{i=1}^{n} (-1)^{i}\left[ T (x_1, x_2, \dots,
 x_{i-1}, x_{i+1},\dots, x_n) \right.
\nonumber \\
&&
- \left. (x_1 \ast x_i, x_2 \ast x_i, \dots, x_{i-1}\ast x_i, x_{i+1}, 
\dots, x_n) \right]
\end{eqnarray}
for $n \geq 2$ 
and $\partial^T_n=0$ for 
$n \leq 1$. 
We regard that the $i=1$ terms contribute $(1-T) (x_2, \ldots, x_n)$. 
Then
 $C_\ast^{\rm TR}(X) = \{C_n^{\rm TR}(X), \partial^T_n \}$
 is a chain complex.
For any $\Lambda$-module $A$, let 
$C_\ast^{\rm TR}(X;A)= 
 \{C_n^{\rm TR}(X) \otimes _{\Lambda} A , \partial^T_n \}$
be the induced chain complex, where the induced boundary operator is 
represented by the same notation.
Let $C^n_{\rm TR}(X;A)=\mbox{Hom}_{\Lambda} ( C_n^{\rm TR}(X), A)$
and define the coboundary operator 
$\delta=\delta^n_{\rm TR}  : C^n_{\rm TR}(X;A) \to C^{n+1}_{\rm TR}(X;A)$
by 
$(\delta f )(c)= (-1)^n f( \partial c) $ 
 for any $c \in  C_n^{\rm TR}(X)$
and $f \in C^n_{\rm TR}(X;A)$. Then 
$C^\ast_{\rm TR}(X;A)= 
 \{C^n_{\rm TR}(X;A) , \delta^n_{\rm TR}  \}$
is a cochain complex.
The $n\/$-th homology and cohomology groups of 
these complexes 
are called 
 {\it twisted rack homology group\/} and {\it cohomology group\/},
and are denoted by $H_n^{\rm TR}(X; A)$ and $H^n_{\rm TR}(X; A)$,
respectively.

\begin{sloppypar}
Let $C_n^{\rm TD}(X;A)$ be the subset of $C_n^{\rm TR}(X;A)$ generated
by $n$-tuples $(x_1, \dots, x_n)$
with $x_{i}=x_{i+1}$ for some $i \in \{1, \dots,n-1\}$ if $n \geq 2$;
otherwise let $C_n^{\rm TD}(X;A)=0$. If $X$ is a quandle, then
$\partial^T_n(C_n^{\rm TD}(X;A)) \subset C_{n-1}^{\rm TD}(X;A)$ and
$C_\ast^{\rm TD}(X;A) = \{ C_n^{\rm TD}(X;A), \partial^T_n \}$ is 
a sub-complex of
$C_\ast^{\rm TR}(X;A)$. 
Similar subcomplexes 
$C^\ast_{\rm TD}(X;A) = \{ C^n_{\rm TD}(X;A), \delta_T^n \}$
are defined for cochain complexes.
\end{sloppypar}

\begin{sloppypar}
The $n\/$-th homology and cohomology groups of 
these complexes  
are called 
 {\it twisted degeneracy  homology group\/} and {\it cohomology group\/},
and are denoted by $H_n^{\rm TD}(X; A)$ and $H^n_{\rm TD}(X; A)$,
respectively.
\end{sloppypar}

Put $C_n^{\rm TQ}(X;A)\! =\! C_n^{\rm TR}(X;A)/ C_n^{\rm TD}(X;A)$ and 
$C_\ast^{\rm TQ}(X;A)\! =\! \{ C_n^{\rm TQ}(X;A), {\partial  }^T_n \}$, 
where all the induced boundary operators are denoted by 
$\partial =  {\partial  }^T_n $. 
A cochain complex 
$C^\ast_{\rm TQ}(X;A) = \{ C^n_{\rm TQ}(X;A), {\delta }_T^n \}$
is similarly defined.
Note again that 
all boundary 
and coboundary 
operators 
will be denoted by $\partial=\partial^T_n$
and $\delta={\delta}_T^n $, respectively. 
The $n\/$-th homology and cohomology groups of 
these complexes 
 are called 
 {\it twisted   homology group\/} and {\it cohomology group\/},
and are denoted by
\begin{eqnarray}
H_n^{\rm TQ}(X; A) 
 = H_{n}(C_\ast^{\rm TQ}(X;A)), \quad
H^n_{\rm TQ}(X; A) 
 = H^{n}(C^\ast_{\rm TQ}(X;A)). \end{eqnarray}
The groups of (co)cycles and (co)boundaries are denoted 
using similar notations.

For ${\rm W}={\rm D},{\rm R},$ or ${\rm Q}$ (denoting the degenerate, 
rack or quandle case, respectively),
the groups of twisted cycles and boundaries are denoted (resp.) by 
${\mbox{\rm ker}}(\partial) =Z_n^{\rm TW}(X;A) \subset C_n^{\rm TW} (X;A)$ and 
${\mbox{\rm Im}}(\partial)=B_n^{\rm TW}(X;A) \subset C_n^{\rm TW}(X;A)$.
 The twisted cocycles and coboundaries are denoted respectively by
${\mbox{\rm ker}}(\delta) =Z^n_{\rm TW}(X;A) \subset C^n_{\rm TW}(X;A)$ and 
${\mbox{\rm Im}}(\partial)=B^n_{\rm TW}(X;A) \subset C^n_{\rm TW}(X;A).$ 
Thus the (co)homology groups are given as quotients:
\begin{eqnarray}
H^{\rm TW}_n(X;A) &=& Z_n^{\rm TW}(X;A)/B_n^{\rm TW}(X;A), \nonumber \\
H_{\rm TW}^n(X;A) &=& Z^n_{\rm TW}(X;A)/B^n_{\rm TW}(X;A). \nonumber
 \end{eqnarray}
See Section~\ref{invsec} for diagrammatic interpretations of
the twisted cycle and cocycle groups.

\begin{example} {\rm
The $1$-cocycle condition is written for $\eta \in Z^1_{ \rm TQ}(X;A)$  
as 
$$ -T \eta(x_2) + T \eta(x_1) + \eta(x_2) - \eta(x_1*x_2)=0, \quad  \mbox{or}$$
$$ T \eta(x_1) + (1-T)\eta(x_2)  = \eta(x_1*x_2) . $$
Note that this means that $\eta: X \rightarrow A$ is a quandle homomorphism.

The $2$-cocycle condition is written for $\phi \in Z^2_{\rm TQ}(X;A)$  as 
\begin{eqnarray*}
\lefteqn{ T [ - \phi (x_2, x_3) + \phi (x_1, x_3) - \phi (x_1, x_2) ] } \\
 & + & [ \phi (x_2, x_3) -  \phi (x_1 * x_2, x_3) 
+ \phi (x_1 * x_3, x_2 * x_3) ] =0 \quad \mbox{or}
\end{eqnarray*}
\begin{eqnarray*}
\lefteqn{ T  \phi (x_1, x_2) +   \phi (x_1 * x_2, x_3) } \\
 & = & T  \phi (x_1, x_3) + (1-T) \phi (x_2, x_3)
 +  \phi (x_1 * x_3, x_2 * x_3).
\end{eqnarray*}
} \end{example}

\begin{example} \label{R3ex} {\rm
We compute $H_2^{\rm TQ}(R_3; R_3)$. 
Let $R_3=\{ 0,1,2\}=\{ a,b,c\}$. 
In this case, note that $R_3=\Z_3 [T, T^{-1}]/(T+1)$,
so $T$ acts as multiplication by $(-1)$, 
 and the boundary homomorphism is computed by
$$ \partial (a,b)= (-1)[ -(b) + (a) ] + [(b)-(a*b)]
= -(a)-(b)-(c) .$$
Since the image is the same for all pair $(a,b)$,
we have $Z_2^{\rm TQ}(R_3; R_3)=(R_3)^5$,
generated by $(a,b)-(0,1)$, $(a,b)\neq (0,1)$.
On the other hand, 
\begin{eqnarray*}
\partial (a,b,a) &=& (-1)[ -(b,a)+(a,a)-(a,b)] \\
 & & - [-(b,a) + (c,a) - (a,c) ] \\
 & = &  (a,b)   -(b,a) + (a,c)   -(c,a), 
\end{eqnarray*}
and 
\begin{eqnarray*}
\partial (a,b,c) &=& (-1) [ -(b,c) + (a,c) - (a,b) ] \\
 & & - [ -(b,c)+ (c,c) - (b,a) ] \\
 & = & (a,b) + (b,a) - (a,c)-(b,c) ,  \\
\end{eqnarray*}
from which it can be seen that 
$\partial (0,1,0)$, $\partial (0,1,2)$, and $\partial (0,2,1)$
span the boundary group   $B_2^{\rm TQ}(R_3; R_3)=(R_3)^3$.
Hence  $H_2^{\rm TQ}(R_3; R_3)=R_3 \times R_3$.
Note that for untwisted case $H_2^{\rm Q}(R_3;A)=0$ for any 
coefficient $A$, see \cite{CJKLS}. 
Also, it can be seen that 
$x= (1,0)-(0,2)$ and $y=(0,1)-(2,1)$
represent generators of  $H_2^{\rm TQ}(R_3; R_3)$.

\begin{sloppypar}
For $A=R_n=\Z_n [T, T^{-1}]/(T+1)$ where $n>3$,
 computations show that
$$ \partial (a,b)= (-1)[ -(b) + (a) ] + [(b)-(a*b)]
= 2(b)-(a)-(c) .$$
 Suppose that $\gcd{(6,n)}=1$. 
Then the boundary map has rank $2$, and
 $Z_2^{\rm TQ}(R_3; R_n)=(R_n)^4$ is 
generated by $e_1= (0,1)+(0,2)+(1,0)$,
$e_2=(0,1)+(0,2)+(2,0)$, 
$e_3=(0,1)-(2,1)$, and
$e_4=(1,2)-(0,2).$
We have
\begin{eqnarray*}
\partial (a,b,a) &=& 2(b,a)  -(c,a) + (a,b) + (a,c)\\
\partial(a,b,c) &=& 2(b,c)-(a,c)+(a,b)+(b,a) . 
\end{eqnarray*}
Substituting various values $\{0,1,2\}$ for $\{a,b,c\}$
 in the above expressions, we obtain:
\begin{eqnarray*}
\partial (0,1,0) &=& 
 2(1,0)  -(2,0) + (0,1) + (0,2)=2 e_1-e_2 \\
\partial (0,2,0) &=& 
 2(2,0)  -(1,0) + (0,2) + (0,1)=2e_2-e_1 \\
\partial(0,1,2) &= & 2(1,2) 
-(0,2)+(0,1) 
+ (1,0) 
=e_1+2e_4\\
\partial(0,2,1) &= & 2(2,1)
-(0,1)+(0,2) 
+ (2,0) 
=e_2-2e_3.
\end{eqnarray*}\end{sloppypar}

Since $\gcd{(n,6)}=1$
and $2$, $3$ and $6$ are invertible in ${\Z}_n$, 
we see that 
$e_1$, $e_2$, $e_3$, and $e_4$ are in the image of
the boundary map. 
Specifically, 
\begin{eqnarray*}
e_1 &=& \partial ( 2/3 
(0,1,0)+1/3(0,2,0))\\
e_2 &=& \partial (1/3(0,1,0)+2/3(0,2,0)) \\
e_3 &=& \partial (1/6(0,1,0)+1/3(0,2,0)-1/2(0,2,1))\\
e_4 &=& \partial (-1/3(0,1,0)+1/2(0,1,2)-1/6(0,2,0)) . 
\end{eqnarray*} 
So $H_2^{\rm TQ}(R_3; R_n)=0$ for any $n>3$, with $\gcd{(n,6)}=1.$
} \end{example}

\begin{example} \label{trivialex} {\rm
Let $X=T_m=\{ 0, 1, \ldots, m-1\}$ be the trivial quandle of $m$ elements, 
so that $a*b=a$ for any $a,b \in X$.
In this case the chain map reduces to 
$(T-1) \partial_0$, where
$$\partial_0 (x_1, \ldots, x_n)
=\sum_{i=1}^n (-1)^i (x_1, \ldots, \widehat{x_i}, \ldots, x_n). $$
In particular, if $T=1$ (in which case the homology is untwisted),
all the chain maps are zero. 
On the other hand, if $(T-1)$ is invertible in the coefficient group $A$,
then the boundary maps coincides with the above $\partial_0$.

For example, we compute $H_2^{\rm TQ}(T_2; A)$ as follows, where assume that
 $(T-1)$ is not a zero divisor.
One computes 
$$ \partial (x,y)=(T-1) [ (x)-(y)] $$
so the kernel $Z_1^{\rm TQ}(T_2; A)$ is written  as 
$\{ a (x,y) + b(y,x) | (T-1) (a-b)=0 \}$. 
Since $(T-1)$ is not a zero divisor, 
this group is the free module generated by $(0,1)+(1,0)$.
On the other hand,
$$ \partial (0,1,0)= (T-1) [ -(1,0)-(0,1) ]
= \partial (1,0,1), $$
so we obtain $H_2^{\rm TQ}(T_2; A)= A / (T-1) A$. 
In particular, if $(T-1)$ is invertible, then $H_2^{\rm TQ}(T_2; A)= 0$.

Cohomology groups are computed similarly,
using characteristic functions.
For example, if $(T-1)$ is not a zero divisor, we find
$H^2_{\rm TQ}(T_2;A)= A / (T-1) A$. 
} \end{example}

The following also follows from the definitions.

\begin{proposition} 
For any quandle $X$ and an Alexander quandle $A$, 
$$H_1^{\rm TQ}(X;A)  \cong A[X]/ ( T x + (1-T) y - x*y ) $$ 
where $A[X]$ is the free module generated by elements of $X$, 
and the quotient is taken by the submodule generated by
elements of the form $ T x + (1-T) y - x*y$ for all $x, y \in X$.
\end{proposition}

\begin{example} {\rm 
For $X=R_3$ and $A=R_n$, $A[X]=R_n (0) \oplus R_n (1) \oplus R_n (2)$,
the free module generated by elements of $R_3$, 
with basis elements denoted by $(0),  (1) $ and $(2)$.
The action by $T$ is multiplication by $(-1)$, 
and the relations $( T x + (1-T) y - x*y ) $
reduce to two of them,
$2(0)-(1)-(2) $ and 
$2(1)-(0)-(2) $. 
These further reduce to $2(0)-(1)-(2) $ and $3[(0)-(1)]$.
Hence $A[X]$ modulo these subgroups is $R_n$ if $(n, 3)=1$,
and $R_n \times R_3$ if $(n,3) \neq 1$. 
Thus we obtain
$$H_1^{\rm TQ}(R_3;R_n) = \left\{ \begin{array}{lr} R_n  & {\mbox{\rm if }} \
(n, 3)=1, \\
R_n \times R_3 & {\mbox{\rm if }} \
(n,3) \neq 1 . \end{array}\right. $$
} \end{example}

\section{Extensions of quandles by Alexander quandles} \label{extsec}

In this section we give interpretations of quandle cocycles 
in low dimensions as 
extensions of quandles. The theories are analogues of   those of group
and other (such as Hochschild) cohomology theories, and are developed 
in parallel to these theories
(see \cite{Brown} Chapter $4$, for example). 

Let $X$ be a quandle and $A$ be an Alexander quandle.
Recall that $\eta\! \in\! Z^1_{\rm TQ}(X;A)$ implies that
$\eta: X \rightarrow A$ is a quandle homomorphism.
Let $0\rightarrow N  \stackrel{i}{\rightarrow} G
  \stackrel{p}{\rightarrow} A \rightarrow 0$ be an exact sequence
of $\Z [T, T^{-1}]$-module homomorphisms among
Alexander quandles. Let $s: A \rightarrow G$ be 
a set-theoretic section
(i.e., $ps=$id$_A$) with the ``normalization condition'' $s(0)=0$.
Then $s\eta: X \rightarrow G$ is a mapping, which is 
not necessarily a quandle homomorphism.
We measure the failure by $2$-cocycles. 
Since $p[T s\eta(x_1)+ (1-T) s\eta (x_2)]=p[ s\eta(x_1*x_2)]$
for any $x_1 , x_2 \in A$, there is $\phi(x_1, x_2) \in N$
such that 
\begin{equation} \label{2cocydefrel}
T s\eta(x_1)+  s\eta (x_2)= i \phi(x_1, x_2) + 
[ T  s\eta (x_2) +  s\eta(x_1*x_2)] . 
\end{equation}
This defines a function $\phi \in C^2_{\rm TQ}(X;N)$.

\begin{lemma} 
$\phi \in Z^2_{\rm TQ}(X;N)$.
\end{lemma}

\proof One computes
\begin{eqnarray*}
\lefteqn{ \underline{T^2 s 
 \eta(x_1) + T s 
 \eta(x_2)} + s\eta(x_3) } \\
& =& [T i \phi(x_1, x_2) + T^2  s\eta (x_2) +  \underline{ T s\eta(x_1*x_2)}]
+  \underline{s  \eta(x_3)} \\
&=& T i \phi(x_1, x_2) + \underline{ T^2  s\eta (x_2)} + 
[i\phi (x_1*x_2, x_3) +   \underline{T s\eta (x_3)} + s\eta( (x_1*x_2)*x_3)] \\
&=& [ T i \phi(x_1, x_2) +i\phi (x_1*x_2, x_3) + T i\phi (x_2, x_3)] \\
 & &  +[ s\eta( (x_1*x_2)*x_3)+ T^2 s \eta(x_3) + T  s\eta (x_2*x_3) ]
\end{eqnarray*}
and on the other hand, 
\begin{eqnarray*}
\lefteqn{ T^2 s 
 \eta(x_1) + \underline{ T s\eta(x_2) + s\eta(x_3)} } \\
& =& [i \phi(x_2, x_3) +  \underline{T  s\eta(x_3)}+s\eta  (x_2*x_3) 
+ \underline{ T^2 s \eta(x_1)} \\
&=& i \phi(x_2, x_3) + [ T i \phi(x_1, x_3) +
T^2 s\eta( x_3) + \underline{ T s \eta(  x_1*x_3)] + s\eta  (x_2*x_3) } \\
&=& [  i \phi(x_2, x_3)+ T i \phi(x_1, x_3) + i \phi( ( x_1*x_3), (x_2*x_3))]\\
& & + [T^2 s\eta( x_3) + T s\eta  (x_2*x_3)+ s \eta(  ( x_1*x_3)*(x_2*x_3))] . 
\quad 
\end{eqnarray*}
The underlines in the equalities indicates where  Relation~(\ref{2cocydefrel})
is going to be applied in the next step of the calculation.
The observant reader will notice that the calculation follows from 
the type III Reidemeister move, compare with Fig.~\ref{2cocy}.
\endproof

Let $s': A \rightarrow G$ be another section, 
and $\phi' \in  Z^2_{\rm TQ}(X;N)$
be a $2$-cocycle determined by 
\begin{equation}
T s'\eta(x_1)+  s'\eta (x_2)= i \phi'(x_1, x_2) + 
[ T  s'\eta (x_2) +  s'\eta(x_1*x_2)] . 
\end{equation}

\begin{lemma}
$[\phi]=[\phi'] \in  H^2_{\rm TQ}(X;N)$.
\end{lemma}
\proof
Since $s'(a)-s(a) \in i(N)$, there is a function $\sigma: A \rightarrow N$
such that $s'(a)=s(a)+i\sigma(a)$ for any $a \in A$. 
Then 
\begin{eqnarray*}
\lefteqn{ T[ s \eta (x_1) + i \sigma \eta(x_1) ] 
+ [ s \eta (x_2) + i \sigma \eta(x_2) ] } \\
&=& i \phi'(x_1, x_2) +
 T [ s \eta (x_2) + i \sigma \eta (x_2) ]
+ [ s \eta (x_1*x_2) + i \sigma \eta (x_1*x_2) ] 
\end{eqnarray*}
and hence $\phi'=\phi - \delta( \sigma \eta ) $.
\endproof

\begin{lemma}
If $[\phi]=0 \in  H^2_{\rm TQ}(X;N)$,
then $\eta: X \rightarrow A $ 
extends to a quandle homomorphism to $G$, 
i.e., there is a quandle homomorphism $\eta': X \rightarrow G$
such that $p \eta'=\eta$.
\end{lemma} 
\proof
By assumption there exists $\xi \in C^1_{\rm TQ}(X;N)$ such that 
$\phi=\delta \xi$. 
By Equality~(\ref{2cocydefrel}), 
the map $\eta'=s \eta - i \xi$ gives rise to a desired
quandle homomorphism.
\endproof

We summarize the above lemmas as follows. 

\begin{theorem} \label{2cocyobstthm}
The obstruction  to extending $\eta: X \rightarrow A$ to 
a quandle homomorphism $X \rightarrow G$ lies
in  $H^2_{\rm TQ}(X;N)$.
\end{theorem}

Such a  $2$-cocycles $\phi$ constructed  above
is called an {\it obstruction $2$-cocycle}.

Next  we use $2$-cocycles to construct extensions.
Let $X$ be a quandle and $A$ be an Alexander quandle.
Let $\phi \in Z^2_{\rm TQ} (X; A)$.
Let $AE(X,A, \phi)$ be the quandle defined on the set $A \times X$ by
the operation 
$(a_1, x_1) * (a_2, x_2) = (a_1 * a_2 + \phi(x_1, x_2), x_1 * x_2)$.

\begin{lemma} \label{prodlemma}
The above defined operation $*$ on $A \times X$ indeed
defines a quandle $AE(X,A, \phi) = (A \times X, *)$,
which is called an {\em Alexander extension} 
of $X$ by $(A, \phi)$. 
\end{lemma}
\proof 
The idempotency is obvious.
For any $(a_2, x_2), (a, x) \in A \times X$,
let $x_1 \in X$ be 
the 
unique element  such that
$x_1 * x_2 = x$
and $a_1 \in A$ be 
the 
unique element   such that
$a_1 * a_2=a - \phi(x_1, x_2)$.
Then it follows that $(a_1, x_1) * (a_2, x_2) =(a, x) $,
and the uniqueness of  $(a_1, x_1)$ with this property is obvious.
The self-distributivity follows from
the $2$-cocycle condition by computation, as follows.
\begin{eqnarray*}
\lefteqn{
   [ (a_1, x_1) * (a_2, x_2)] * (a_3, x_3) } \\
  &=&
( Ta_1 + (1-T)a_2  +  \phi(x_1, x_2), x_1 * x_2) * (a_3, x_3) \\
&=& (  T[ Ta_1 + (1-T)a_2 + \phi(x_1, x_2)] \\ & & + (1-T)a_3 +
  \phi(x_1 * x_2,  x_3) ,  (x_1 * x_2)* x_3 ) \\
&=& ( ( a_1*a_2)*a_3 + T  \phi(x_1, x_2)\\ & & +  \phi(x_1 * x_2,  x_3),
 (x_1 * x_2)* x_3 ), 
\end{eqnarray*}
and
\begin{eqnarray*}
\lefteqn{
  [ (a_1, x_1) * (a_3, x_3)] *[ (a_2, x_2) * (a_3, x_3)]} \\
  &=&
  (a_1*a_3 +  \phi(x_1, x_3), x_1 * x_3) *  (a_2*a_3 + 
 \phi(x_2, x_3), x_2 * x_3) \\
&=& ( T [ a_1*a_3 +  \phi(x_1, x_3)] +
 (1-T) [ a_2*a_3 + \phi(x_2, x_3)] \\ & & + \phi( x_1 * x_3,   x_2 * x_3) ,
(x_1 * x_3) *(x_2 * x_3) ) \\
&=& 
( (a_1*a_3)*(a_2*a_3) + T \phi(x_1, x_3) +  (1-T) \phi(x_2, x_3) \\& & + 
 \phi( x_1 * x_3,   x_2 * x_3) ,(x_1 * x_3) *(x_2 * x_3) ).
\end{eqnarray*}
They are equal by the $2$-cocycle condition.
\endproof

\begin{remark} \label{obstrem} {\rm
In Theorem~\ref{2cocyobstthm}, we consider the situation where $A=X$,
$\eta=$id, and $G=E=AE(X,B,\phi)$  
 for some cocycle $\phi \in Z^2_{\rm TQ}(X;B)$
where $X$, $E$, $B$ are Alexander quandles.

Assume that we have a short exact sequence 
$$ 0 \rightarrow B \stackrel{i}{\rightarrow} E 
\stackrel{p}{\rightarrow} X \rightarrow 0 $$
of $\Z[T, T^{-1}]$-modules, where $i(b)=(b, 0)$ and $p((b,x))=x$
for $b \in B$ and $(b, x) \in E=B \times X$. 
Then there is a 
section $s: X \rightarrow E$ defined by $s(x)=(0, x)$
satisfying $ps = $id. 
Then we have 
\begin{eqnarray*}
\lefteqn{  [ T s (x_1) + (1-T) s(x_2) ] - s(x_1 * x_2) } \\
 & = & s(x_1) * s(x_2) -  s(x_1 * x_2) \\
 &=& (0, x_1)*(0, x_2) - (0, x_1 * x_2) \\
 &=& ( \phi(x_1, x_2), 0) \\
&=& i\phi(x_1, x_2).
\end{eqnarray*}
Therefore the cocycle used in the preceding Lemma,
which we may call an {\it extension cocycle},  is an obstruction cocycle.
} \end{remark}

\begin{definition}{\rm
Two surjective homomorphisms of quandles $\pi_j: E_j \rightarrow X$,
$j=1,2$,
are called {\it equivalent } if there is
a quandle isomorphism $f: E_1 \rightarrow E_2$
such that $\pi_1=\pi_2 f$.}
  \end{definition}

Note that there is a natural surjective homomorphism
$\pi: AE(X,A, \phi)=A \times X  \rightarrow X$, which is the projection to
the second factor.

\begin{lemma} \label{equivlemma}
If   $\phi_1$ and $\phi_2$ are cohomologous, i.e.,
$[\phi_1]=[\phi_2] \in  H^2_{\rm TQ}(X;A)$,
then
$\pi_1: AE(X,A, \phi_1) \rightarrow X$ and  $\pi_2: AE(X,A, \phi_2)\rightarrow X$ 
are equivalent.
\end{lemma}
\proof
There is a $1$-cochain $\eta \in C^1_{\rm TQ}(X;A)$ such that
$\phi_1 =\phi_2 + \delta \eta$.
We show that
  $f: AE(X,A, \phi_1)=A \times X \rightarrow A \times X = AE(X,A, \phi_2)$
defined by $f(a,x)=(a + \eta(x), x) $ gives rise to an equivalence.
First we compute
\begin{eqnarray*}
  f( (a_1, x_1) * (a_2, x_2) ) &  = &
 f( (a_1* a_2 + \phi_1( x_1, x_2 ) ,  x_1 *  x_2 ) ) \\
  & = &  (a_1* a_2 + \phi_1( x_1, x_2 ) + \eta( x_1 * x_2 ) ,  x_1 * x_2 ), \; \; 
\mbox{and} \\
  f( (a_1, x_1) ) * f(  (a_2, x_2) )  & = & (a_1 + \eta(x_1) , x_1)
 * (a_2  + \eta (x_2) ,  x_2) \\
&=&  ( T( a_1 + \eta(x_1)) + (1-T) (a_2 + \eta(x_2)) \\ & & + \phi_2(x_1, x_2) , 
  x_1 * x_2 ) \\
& = & ( a_1*a_2 +  \phi_2(x_1, x_2) \\ & &  + ( T \eta(x_1) +  (1-T)\eta(x_2)),
  x_1 * x_2 )
\end{eqnarray*}
which are equal since $\phi_1 =\phi_2 + \delta \eta$.
Hence $f$ defines a quandle homomorphism.
The map $f': A \times X \rightarrow A \times X$ defined by
$f'(a, x)=(a - \eta(x), x)$ defines the inverse of $f$,
hence $f$ is an isomorphism. The map $f$ satisfies
$\pi_1=\pi_2 f$ by definition.
\endproof

\begin{lemma}
If
natural surjective homomorphisms 
(the projections to the second factor
$A \times X \rightarrow X$)      
$AE(X,A, \phi_1)\rightarrow X$ and  $AE(X,A, \phi_2)\rightarrow X$
  are equivalent, then  $\phi_1$ and $\phi_2$ are cohomologous:
$[\phi_1]=[\phi_2] \in  H^2_{\rm TQ}(X;A)$.
\end{lemma}
\proof 
Let $f:  AE(X,A, \phi_1)=A \times X \rightarrow A \times X = AE(X,A, \phi_2)$
be a quandle isomorphism with $\pi_1=\pi_2 f$.
Since $\pi_1(a,x)=x= \pi_2( f(a, x))$,
  there is an element $\eta(x) \in A$ such that
$f(a, x)= (a + \eta(x), x)$,  for any $x \in X$.
This defines a function $\eta: X \rightarrow A $, $\eta \in C^1_{\rm TQ}(X;A)$.
The condition that $f$ is a quandle homomorphism
implies that  $\phi_1 =\phi_2 + \delta \eta$ by
the same computation as the preceding lemma.
Hence the result follows.
\endproof

The lemmas imply the following theorem.

\begin{theorem} \label{2cocythm}
There is a bijection between the equivalence classes of
natural surjective homomorphisms 
  $AE(X,A, \phi)\rightarrow X$ for a fixed $X$ and $A$,
   and
the set $H^2_{\rm TQ}(X;A)$.
\end{theorem}

Next we consider interpretations of $3$-cycles in extensions of quandles.
Let
$0 \rightarrow N   \stackrel{i}{\rightarrow} G  \stackrel{p}{\rightarrow} A 
\rightarrow 0$
be a short exact sequence of $\Z [T,T^{-1}]$-modules.
Let $X$ be a quandle. For $\phi \in Z^2_{\rm TQ}(X;A)$,
let $AE(X,A, \phi)$ be as above.
Let $s: A \rightarrow G$ be a set-theoretic
(not necessarily group homomorphism)
section,
i.e., $ps=\mbox{id}_A$,
with the ``normalization condition''  of $s(0)=0$.  

Consider the binary operation
$ (G \times X ) \times ( G \times X)
\rightarrow G \times X$ defined by
\begin{equation} \label{extdef}
  (g_1, x_1) * (g_2, x_2) = (g_1* g_2 +  s \phi (x_1, x_2) , x_1 * x_2 ).
\end{equation}
  We describe an obstruction to this being a quandle operation
by $3$-cocycles.

Since $\phi $ satisfies the $2$-cocycle condition,
$$p (T  s \phi(x_1, x_2) +  s \phi (x_1 * x_2, x_3)) $$ 
$$= p ( T  s \phi (x_1, x_3) +  (1-T) s \phi ( x_2, x_3) +
s \phi (x_1 * x_3, x_2 * x_3 ) ) $$
  in $A$.
Hence there is a function
  $\theta: X \times X \times X \rightarrow N$ such that
\begin{equation} \label{3cocydef}
\begin{array}{l}
  T s \phi (x_1, x_2) + s \phi (x_1 * x_2, x_3) + T  s \phi ( x_2, x_3) \\
=i \theta (x_1, x_2, x_3) +  s \phi ( x_2, x_3)+ T s \phi (x_1, x_3) 
 + s \phi (x_1 * x_3, x_2 * x_3 ),
\end{array}
\end{equation}
where we moved the term  $T  s \phi ( x_2, x_3)$
so that we have only positive terms in the definition of $\theta$.

\begin{lemma} \label{movielemma}
$\theta \in Z^3_{\rm Q}(X;N)$.
\end{lemma}
\proof 
First, if $x_1=x_2$, or $x_2=x_3$, then the above defining relation 
for $\theta$ implies that $\theta (x_1, x_1, x_3)=1=\theta(x_1, x_2, x_2)$. 
For the $3$-cocycle condition, 
one computes 
\begin{eqnarray*}
\lefteqn{ \underline{ T^2 s \phi (x_1, x_2) + T s \phi (x_1 * x_2, x_3)
+ T^2 s \phi( x_2, x_3) } } \\ & &
 + s\phi ( (x_1 * x_2)*x_3, x_4)
+ T s\phi (x_2 * x_3 ,  x_4) + T^2 s\phi( x_3, x_4)  \\
 &=& i T \theta (x_1, x_2, x_3) \\ && + [  T s \phi( x_2, x_3)
+ T^2 s \phi (x_1, x_3) + \underline{ Ts \phi (x_1 * x_3, x_2 * x_3 ) } 
] \\ & &
 +  \underline{s \phi ( (x_1 * x_2)*x_3, x_4) 
+  T s\phi (x_2 * x_3 ,  x_4)} +  T^2 s\phi( x_3, x_4)  
\end{eqnarray*}
\begin{eqnarray*}
 &=& [i T \theta (x_1, x_2, x_3) + i \theta (x_1 *x_3, x_2 * x_3, x_4) ] \\ & &
 + [ \underline{ s \phi(x_2 * x_3, x_4) + T  s \phi (x_1 * x_3,  x_4)} 
\\ & &+  s\phi ((x_1 *x_3)*x_4,  (x_2 * x_3  )*x_4 )] \\ & &
+  \underline{ T^2 s \phi (x_1, x_3)} + T  s \phi( x_2, x_3)
 +  \underline{ T^2  s\phi( x_3, x_4)}  \\
 &=& [i T \theta (x_1, x_2, x_3) + i \theta (x_1 *x_3, x_2 * x_3, x_4)
+ i T \theta(x_1, x_3, x_4)] \\ & &
 + [ \underline{Ts \phi(x_3, x_4)} + T^2 s\phi (x_1,  x_4) + T s\phi (x_1*x_4,  x_3*x_4) ]
\\ & &
+  s\phi ((x_1 *x_3)*x_4,  (x_2 * x_3 )*x_4 )
+  \underline{ T s \phi( x_2, x_3)
+  s\phi (x_2 * x_3 ,  x_4)  }  \\
 &=& [i T \theta (x_1, x_2, x_3) +  i\theta (x_1 *x_3, x_2 * x_3, x_4) \\ & & 
+ i T \theta(x_1, x_3, x_4) + i  \theta( x_2, x_3, x_4) ] \\ & &
+ [s \phi(x_3, x_4) + T s \phi( x_2, x_4) + s \phi( x_2*x_4, x_3*x_4) ] \\ & &
+ T^2 s\phi (x_1,  x_4) + T s\phi (x_1*x_4,  x_3*x_4) \\ & &
+ s\phi ((x_1 *x_3)*x_4,  (x_2 * x_3 )*x_4 )
\end{eqnarray*}
and on the other hand,
\begin{eqnarray*}
\lefteqn{ T^2 s \phi (x_1, x_2) + T  s \phi (x_1 * x_2, x_3)
 + \underline{T^2 s \phi( x_2, x_3)} }  \\ & & 
 +  s\phi ( (x_1 * x_2)*x_3, x_4)
+ \underline{ T s\phi (x_2 * x_3 ,  x_4)  + T^2  s\phi( x_3, x_4) }   \\
 &=& i T \theta( x_2, x_3, x_4) + [ \underline{ Ts \phi(x_3, x_4)} + 
 T^2  s \phi( x_2, x_4)+T  s\phi (x_2 * x_4 , x_3* x_4)] \\ & &
+ T^2  s \phi (x_1, x_2) + \underline{ T s \phi (x_1 * x_2, x_3)
    + s\phi ( (x_1 * x_2)*x_3, x_4) } \\
  &=&  [i T \theta ( x_2, x_3, x_4) + i\theta (x_1 * x_2, x_3, x_4) ]
 \\ & & 
+ [ s \phi ( x_3, x_4) + 
\underline{ T s \phi (x_1 * x_2, x_4) } + s\phi ( (x_1 * x_2)*x_4, x_3*x_4 ) ]
\\ & &
+ \underline{T^2  s \phi( x_2, x_4)} 
+ T s\phi (x_2 * x_4 , x_3* x_4) + \underline{T^2  s \phi (x_1, x_2)} \\
 &=& [i T \theta( x_2, x_3, x_4) +  i\theta(x_1 * x_2, x_3, x_4)
+ i T \theta( x_1, x_2, x_4) ]  \\ & &
+ [ Ts \phi(x_2, x_4) 
+  T^2 s \phi (x_1, x_4) +  \underline{ T s \phi (x_1 * x_4 ,  x_2* x_4) } ] \\ & &
+  \underline{ s\phi ( (x_1 * x_2)*x_4, x_3*x_4 ) }
 + \underline{ T s\phi (x_2 * x_4 , x_3* x_4) } 
+   s\phi( x_3, x_4)  \\
   &=&  [i T \theta( x_2, x_3, x_4) + i\theta(x_1 * x_2, x_3, x_4) \\ & &
+ iT \theta( x_1, x_2, x_4) + i\theta(x_1 * x_4, x_2 * x_4, x_3 * x_4) ]
\\ & &
+ [ s \phi(x_2*x_4, x_3*x_4) 
+  T s \phi (x_1 * x_4 ,  x_3 * x_4)
\\ && +  s \phi ( (x_1 * x_3)*x_4,  (x_2 * x_3) * x_4)] \\ & &
+  T^2 s \phi (x_1, x_4) + Ts \phi( x_2, x_4) 
+   s\phi( x_3, x_4)
\end{eqnarray*}
so that we obtain the result.
The underlines in the equalities indicate where the relation (\ref{3cocydef})
is going to be applied in the next step of the calculation.
\endproof

The above computation was facilitated by knot diagrams colored by 
quandle elements, and their movies, by a direct correspondence.
This diagrammatic  
method of computations is discussed in 
Section~\ref{invsec}.

Let $s': A \rightarrow G$ be another section, and  $\theta ' $
be a $3$-cocycle defined similarly for $s'$ by
\begin{equation} \label{3cocydef2}
\begin{array}{l}
 T s' \phi (x_1, x_2) +   s' \phi (x_1 * x_2, x_3) + T s' \phi (x_2, x_3) \\
=i \theta ' (x_1, x_2, x_3) + T s' \phi (x_1, x_3)
+  s' \phi (x_2, x_3) +  s' \phi (x_1 * x_3, x_2 * x_3 ).
\end{array}
\end{equation}

\begin{lemma}
The two $3$-cocycles $\theta$ and $\theta'$ are cohomologous,
$[\theta]=[\theta'] \in H^3_{\rm TQ}(X;N)$.
\end{lemma}
\proof
Since $ s'(a) - s(a)  \in i (N) $
for any $a \in A$,
there is a function $\sigma: A \rightarrow N$ such that
$s'(a) =s(a) +  i \sigma(a)$ for any $a \in A$.
{}From Equality~(\ref{3cocydef2}) we obtain
\begin{eqnarray*}
\lefteqn{
  T [ s \phi (x_1, x_2) + i \sigma \phi (x_1, x_2) ] 
 +  [s \phi (x_1 * x_2, x_3) +  i \sigma\phi (x_1 * x_2, x_3) ] } \\ & & 
+ T[ s \phi (x_2, x_3) + i \sigma \phi  (x_2, x_3)]  \\
&=&
i \theta ' (x_1, x_2, x_3)
 + T [ s \phi (x_1, x_3) + i \sigma \phi (x_1, x_3) ] \\
 & & + [ s \phi (x_2, x_3) + i \sigma \phi  (x_2, x_3)] \\
& & 
 + [s \phi (x_1 * x_3, x_2 * x_3 ) + i \sigma \phi (x_1 * x_3, x_2 * x_3 )].
\end{eqnarray*}
Hence we have $\theta'= \theta + \delta (\sigma \phi)$.
\endproof 

\begin{lemma}
If $\theta$ is a coboundary, i.e.,
$[\theta ] =0 \in  H^3_{\rm TQ}(X;N)$, then
$G \times X$ admits a quandle structure
such that $p \times \mbox{id}_X: G \times X \rightarrow A \times X$
is a quandle homomorphism.
\end{lemma}
\proof 
By assumption there is $\xi \in C^2_{\rm TQ}(X;N)$ such that
$\theta=\delta \xi$.
Define a binary operation on $G \times X$ by
$$
(g_1, x_1)*(g_2, x_2) = (g_1 *g_2 +   s \phi(x_1, x_2) - i \xi (x_1, x_2),
x_1 * x_2 ). $$
Then by Equality~(\ref{3cocydef}),
this defines a desired quandle operation.
\endproof

We summarize the above lemmas as

\begin{theorem} \label{obstthm}
The obstruction to extending  the quandle $AE(X,A, \phi)=A \times X$ to
$G \times X$ lies in $H^3_{\rm TQ}(X;N)$.
  \end{theorem}

Such a  $3$-cocycles $\theta$ constructed  above
is called an {\it obstruction $3$-cocycle}.

\section{Alexander quandles as Alexander extensions} \label{lxsec}

\begin{lemma} \label{cocylemma}
Let $X$, $E$ be  quandles, and $A$ be an Alexander quandle.
Suppose there exists a bijection 
$f: E \rightarrow A \times X$ with the following property.
There exists a function $\phi: X \times X \rightarrow A$ such that
for any $e_i \in E$ ($i=1,2$), 
if $f(e_i)=(a_i, x_i)$, then 
$f(e_1 *  e_2) = (a_1 *a_2 + \phi(x_1, x_2) , x_1 * x_2 )$. 
Then $\phi \in Z^2_{\rm TQ}(X; A)$. 
\end{lemma} 
\proof 
For any $x \in X$ and $a \in A$,
 there is $e \in E$ such that   $f(e)=(a,x)$, and 
$$ (a,x)=f(e)=f(e*e)=(a *a + \phi(x, x), x) , $$
so that we have $\phi(x, x)=0$ for any $x \in X$.

By identifying $A \times X$ with $E$ by $f$, the quandle operation
$*$ on  $A \times X$ is defined, for any  $(a_i, x_i)$ ($i=1,2$), by
 $$(a_1, x_1) * (a_2, x_2) = (a_1*a_2 +  \phi(x_1, x_2), x_1 * x_2).$$
Since  $A \times X$ is a quandle isomorphic to $E$ under this $*$, we have 
\begin{eqnarray*}
\lefteqn{
  [ (a_1, x_1) * (a_2, x_2)] * (a_3, x_3) } \\
 &=&
(a_1 * a_2 +  \phi(x_1, x_2), x_1 * x_2) * (a_3, x_3) \\
&=& ( (a_1*a_2) *a_3 + T \phi(x_1, x_2) + \phi(x_1 * x_2,  x_3) ,  (x_1 * x_2)* x_3 ) , 
\end{eqnarray*} 
and
\begin{eqnarray*}
\lefteqn{
 [ (a_1, x_1) * (a_3, x_3)] *[ (a_2, x_2) * (a_3, x_3)]} \\
 &=& 
 (a_1*a_3 +  \phi(x_1, x_3), x_1 * x_3) *  (a_2*a_3 +  \phi(x_2, x_3), x_2 * x_3) \\
&=& (  (a_1*a_3)*(a_2*a_3 ) + T \phi(x_1, x_3) 
\\ &  & 
+ (1-T) \phi( x_2, x_3) + \phi( x_1 * x_3,   x_2 * x_3), 
(x_1 * x_3) *(x_2 * x_3) )
\end{eqnarray*}
are equal for any  $(a_i, x_i)$ ($i=1,2,3$).
Hence $\phi$ satisfies the $2$-cocycle condition.
\endproof

This lemma implies that under the same assumption  
we have $E=AE(X,A,\phi)$,  where  $\phi \in Z^2_{\rm TQ}(X; A)$.
Next we identify such examples.

Let $\Lambda_p=\Z_p[T, T^{-1}]$ for a positive integer $p$
(or $p=0$, in which case  $\Lambda_p$ is understood to be 
$\Lambda=\Z[T, T^{-1}]$). 
Note that  since $T$ is a unit in $\Lambda_p$,
$\Lambda_p/(h)$ for a Laurent polynomial $h \in \Lambda_p$ is
isomorphic to $\Lambda_p/(T^n h)$ for any integer $n$, so that 
we may assume that $h$ is a polynomial with a non-zero constant 
(without negative exponents of $T$).

\begin{lemma} \label{coeffcocylemma}
Let $h \in \Lambda_{p^{m}}$ be a polynomial with the leading and 
constant coefficients invertible, or $h=0$.
 Let $\bar{h} \in  \Lambda_{p^{m-1}}$
and $\tilde{h} \in  \Lambda_{p}$ be such that 
$\bar{h} \equiv h \ \mbox{\rm mod} \ (p^{m-1})$ 
and $\tilde{h} \equiv h \ \mbox{\rm mod} \ (p)$, respectively
(in other words, $\bar{h}$ is $h$ with its
coefficients reduced modulo $p^{m-1}$, and $\tilde{h}$ is $h$
 with its coefficients reduced modulo $p$).
Then the   quandle $E=\Lambda_{p^m}/({h})$
satisfies the conditions in Lemma~\ref{cocylemma}
with $X=\Lambda_{p^{m-1}}/(\bar{h})$ and $A=\Lambda_p / (\tilde{h})$. 

In particular, $\Lambda_{p^m}/({h}) $ is an Alexander extension
of $\Lambda_{p^{m-1}}/(\bar{h})$ by
 $\Lambda_p / (\tilde{h})$:
$$\Lambda_{p^m}/(h)
 = AE(\Lambda_{p^{m-1}}/(\bar{h}), \ \Lambda_p / (\tilde{h}) , \ \phi) , $$
for some 
$\phi \in Z^2_{\rm TQ}(\Lambda_{p^{m-1}}/(\bar{h}); \Lambda_p / (\tilde{h}) )$.
\end{lemma}
\proof 
Let $A \in \Z _{p^m}$. Represent $A$ in $p^m$-ary notation
as $$A=\sum_{i=0}^{m-1} A_i p^i$$ where 
 $A_i \in \{0, \ldots, p-1\}.$
Since $p$ is fixed throughout, we  represent $A$ by the sequence
$$[A_{m-1}, A_{m-2}, A_{m-3}, \ldots, A_0 ]. $$ 
Define $\overline{A} = [A_{m-2}, \ldots, A_0].$
Observe that $A \equiv \overline{A} \pmod{ p^{m-1}}$, and
$A \equiv A_0 \pmod{p}$. 

 Let  $\hat{\pi}: \Z_{p^m} \rightarrow \Z_{p^{m-1}}$ be the map 
defined by $\hat{\pi}(A) = \overline{A}$. We obtain a short exact
sequence:
$$0 \rightarrow \Z_p \stackrel{\hat{\imath}}{\rightarrow} \Z_{p^m}
\stackrel{\hat{\pi}}{\rightarrow} \Z_{p^{m-1}} \rightarrow 0$$
where $\hat{\imath}(A)=[A,0,\dots,0]$. There is a set-theoretic section
$ \Z_{p^m} \stackrel{\hat{s}}{\leftarrow} \Z_{p^{m-1}}$ defined by
$\hat{s}[A_{m-2},\dots,A_0]=[0,A_{m-2},\dots,A_0].$
The map $\hat{s}$ satisfies $\hat{\pi}\hat{s}={\rm id}$ and 
$\hat{s}(0)=0$.

For a polynomial $L(T) \in \Lambda_{p^m} = \Z_{p^m}[T, T^{-1}]$,
write $$L(T) = \sum_{j=-n}^k [A_{j,m-1}, A_{j, m-2}, \ldots, A_{j,0}]
T^j.$$
Define $$\overline{L}(T) = \sum_{j=-n}^k [ A_{j, m-2}, \ldots, A_{j,0}] T^j \in
\Lambda_{p^{m-1}},$$ 
 and 
$$\tilde{L}(T) = \sum_{j=-n}^k A_{j,m-1} T^j \in \Lambda_{p}.$$
There is a one-to-one correspondence $f: \Lambda_{p^m} \rightarrow
\Lambda_{p} \times \Lambda_{p^{m-1}}$ given by
$f(L)= (\tilde{L}, \overline{L})$. 
We have a short exact sequence of rings:
$$0 \rightarrow \Z_p[T,T^{-1}] \stackrel{i}{\rightarrow}
\Z_{p^m}[T,T^{-1}] \stackrel{\pi}{\rightarrow} \Z_{p^{m-1}}[T,T^{-1}]
\rightarrow 0$$
with a set theoretic section
$ \Z_{p^m}[T,T^{-1}] \stackrel{s}{\leftarrow} \Z_{p^{m-1}}[T,T^{-1}]$
where $i$, $\pi$ and $s$ are the natural maps induced by 
$\hat{i}$, $\hat{\pi}$ and $\hat{s}$, respectively.
Note that 
for $L \in \Lambda_{p^m} = \Z_{p^m}[T,T^{-1}]$ we have
 $\overline{L}= \pi(L)$, and 
the section $s:\Lambda_{p^{m-1}} \rightarrow \Lambda_{p^m}$ is defined by
the formula
$$s \left(\sum_{j=-n}^k [ A_{j, m-2}, \ldots, A_{j,0}] T^j\right) = 
\sum_{j=-n}^k [0, A_{j, m-2}, \ldots, A_{j,0}] T^j.$$
For $L, M \in \Lambda_{p^m}$, let  
$$s(\overline{L})*s(\overline{M})= \sum_{j} [F_{j, m-1}, \ldots, F_{j, 0}] T^j \in \Lambda_{p^{m-1}}.$$
If $L= \sum_j A_j T^j$, and $M= \sum_j B_j T^j$, then 
 $$L*M =  B_{-n} T^{-n} + \sum_{j=-n+1}^{k+1} (A_{j-1} - B_{j-1} +B_j)
T^j = \sum_{j=-n}^k C_j T^j.$$
Furthermore, 
\begin{eqnarray*}
\overline{L} *\overline{M}  &=& [B_{-n,m-2}, \ldots, B_{-n,0}] T^{-n} \\
           &+& \sum_{j=-n+1}^{k+1} \left( [A_{j-1,m-2}, \ldots ,
A_{j-1,0}] \right. \\ & -& \left. [B_{j-1,m-2}, \ldots , B_{j-1,0}] +[B_{j,m-2}, \ldots ,
B_{j,0}] \right) T^j
\end{eqnarray*} 
and write the right-hand side by $ \sum_{j=-n}^k D_j T^j$.
Note that $D_j$'s are well-defined integers, 
not only  elements of $\Z_{p^{m-2}}$.
If $D_j$ is positive, then $F_{j, m-1}=0$, and if 
 $D_j$ is negative, then $F_{j, m-1}=p-1$.
Hence
$$f(L*M)=(\tilde{L}*\tilde{M} +
\phi(\overline{L} , \overline{M}), \overline{L}*\overline{M}), $$
where 
$$ \phi( \overline{L}, \overline{M})=\sum_{j=-n}^k F_{j, m-1}. $$
This concludes the case $h=0$.

Now let $h(T) \in \Z_{p^m}[T]$  be a  polynomial with the leading and constant
coefficients being 
invertible in $\Z_p$. 
  Let $(h)$ denote the ideal generated by $h$.
Since  $i(\tilde{h}) \subset (h)$,  
we obtain a short exact
sequence of quotients:
$$0 \rightarrow \Z_p[T,T^{-1}]/(\tilde{h}) \stackrel{\overline{\imath}}{\rightarrow}
\Z_{p^m}[T,T^{-1}]/(h) \stackrel{\overline{\pi}}{\rightarrow}
\Z_{p^{m-1}}[T,T^{-1}]/(\overline{h}) \rightarrow 0$$
with a set-theoretic section 
$ \Z_{p^m}[T,T^{-1}]/(h) \stackrel{\overline{s}}{\leftarrow}
\Z_{p^{m-1}}[T,T^{-1}]/(\overline{h}).$
Thus we obtain a twisted cocycle 
$$\phi: \Z_{p^{m-1}}[T,T^{-1}]/(\overline{h}) \times
\Z_{p^{m-1}}[T,T^{-1}]/(\overline{h}) \rightarrow \Z_p[T,T^{-1}]/(\tilde{h}).
\eqno{\qed} $$

Since $R_n=\Lambda_n / (T+1)$, we have the following.

\begin{corollary} \label{dihedcocylemma}
The dihedral quandle $E=R_{p^m}$, where $p,m$ are positive integers with $m>1$,
satisfies the conditions in Lemma~\ref{cocylemma}
with $X=R_{p^{m-1}}$ and $A=R_p$. 

In particular, $R_{p^m} $ is an Alexander extension
of $R_{p^{m-1}}$ by $R_p$: $$R_{p^m}= AE(R_{p^{m-1}},   R_p,  \phi),$$ 
for some $\phi \in Z^2_{\rm TQ}(R_{p^{m-1}}; R_p)$.
\end{corollary}

\begin{example} \label{dihedex} {\rm
Let $X=R_3$ and $A=R_3$, then the proof of Lemma~\ref{coeffcocylemma}
gives an explicit $2$-cocycle $\phi$ as follows. For 
$\phi(r_1, r_2)=\phi(1,2)$, for example, one computes
$$ r_1 * r_2 = [0,1]* [0,2]= 2 [0,2]- [0,1] =3 =3 \cdot 1 + 0=[1,0],$$ 
 Hence $\phi(0,2)=1$. 
In terms of  the characteristic function, 
the cocycle $\phi$ contains the term $\chi_{0,2}$, where 
$$\chi_{a,b} (x,y) = \left\{ \begin{array}{ll} 1 & {\mbox{\rm if }} \
(x,y)=(a,b), \\
0 & {\mbox{\rm if }} \
(x,y)\not=(a,b) \end{array}\right.$$
is the characteristic function.
By computing the quotients for all pairs, one obtains
$$ \phi = \chi_{0,2} + \chi_{1,2} +  2  \chi_{1,0} + 2 \chi_{2,0}. $$
} \end{example}

\begin{proposition} \label{integer}
The quandle $R_{\infty} $ is an Alexander extension of $R_n$ by 
 $R_{\infty} $, for any positive integer $n$.
\end{proposition}
\proof
Consider the short exact sequence of abelian groups:
$$0 \rightarrow \Z \stackrel{ \cdot n }{\rightarrow} \Z
\stackrel{\pi}{\rightarrow} \Z_n \rightarrow 0 $$
The groups $\Z$ and $\Z_n$ are quandles under the operation:
$a*b= 2b -a$. In the latter case the quantity $2b-a$ is interpreted modulo
$n$. In the former case, it is an integer. The quandle $R_n$ is the
set $\Z_n=\{0, \ldots, n-1 \}$ with this operation. We can define
 a set-theoretic section $s: R_n \rightarrow \Z$ by $s(a)=a$. 
For $a \in \Z$, let $a=\tilde{a} n + \overline{a}$, where
$\tilde{a} \in \Z$ and $0 \leq \overline{a} < n$ are the quotient
and remainder.
Define $f:  \Z \rightarrow E=\Z \times \Z_n $ by 
$f(a) = ( \tilde{a}, \;  \overline{a} \  \mbox{mod} \ (n) )$. 
Write  $s(a)*s(b)=2  \overline{b}- \overline{a} = q n + r$
where $q \in \Z$ and $0 \leq r < n$. 
Then 
$$f(a*b) = f(2b-a) = (2 \tilde{b} - \tilde{a} )n +  (q n + r)
 =(2 \tilde{b} - \tilde{a} +q )n + r, $$
so that we have 
$$f(a*b)=(\tilde{a} * \tilde{b} + \phi (\overline{a},  \overline{b}),
\overline{a} *  \overline{b} \  \mbox{mod} \ (n) ) . $$
The cocycle $\phi$ is given by
$$\phi(a,b)= \left\{ \begin{array}{rrc} 
-1 & {\mbox{\rm if}} & 2b < a, \\
 0 & {\mbox{\rm if}} & 2b < n+a  \ \ {\mbox{\rm and}} \ \ a \le 2b, \\
1 & {\mbox{\rm if}} & n+a \le 2b. \end{array} \right.      $$
Thus in terms of characteristic functions:
$$\phi =  \sum_{n+a \le 2b} \chi_{a,b} - \sum_{2b<a} \chi_{a,b} \eqno{\qed}$$

\begin{example} {\rm
For $R_3$, we obtain 
$$ \phi = \chi_{0,2} + \chi_{1,2} -  \chi_{1,0} -  \chi_{2,0}. $$
} \end{example}

\begin{proposition} The cocycle $\phi \in Z^2_{\rm Q} (R_n;R_\infty)$
given in Proposition~\ref{integer} is not a coboundary. \end{proposition}
 \proof By Lemma~\ref{equivlemma}, if $\phi$ were a
coboundary, then $R_{\infty} $  
would be 
isomorphic to 
$R_{\infty} \times R_n$, which contains a finite subquandle $R_n$. 
A finite subquandle of $R_\infty$ has a largest element $M$. Let $a$ be
any
other element; then $2M-a \le M$, so
$a=M$.
Hence the only finite subquandles of
$R_\infty$ are the $1$-element trivial quandles.  \endproof

\begin{theorem}
Let $h \in \Lambda_{n}$ be a polynomial with the leading and constant 
coefficients invertible.
Let $\Lambda_{n}/(h)$ be a dihedral quandle, where $n$ is a positive integer
with the prime decomposition $n=p_1^{e_1} \ldots p_k^{e_k}$,
for a positive integers $e_1, \ldots, e_k$ and $k$. 

Then as quandles $ \Lambda_n / (h)$ is isomorphic to 
$ \Lambda_{p_1^{e_1}}/(h_1) \times \ldots \times  \Lambda_{p_k^{e_k}}/(h_k) $,
where $h \equiv h_j \ \mbox{\rm mod} \ (p_j^{e_j})$, 
and each factor $ \Lambda_{p_j^{e_j}} / (h_j)$ is inductively described
as 
an 
Alexander extension: $$ \Lambda_{p_j^{d_j}}/(h_j)
 = AE(  \Lambda_{p_j^{d_j -1}} / (\bar{h_j}),\   \Lambda_{p} / (\tilde{h_j}), \ \phi) , $$
for some $\phi \in Z^2_{\rm TQ}( \Lambda_{p_j^{d_j -1}} / (\tilde{h_j}) ;  \Lambda_{p} / (\tilde{h_j}))$,
where $h_j \equiv \bar{h_j} \ \mbox{\rm mod} \ (p_j^{d_j-1})$
and $h_j \equiv \bar{h_j}'  \ \mbox{\rm mod} \ (p_j)$.
\end{theorem}
\proof
As 
rings, 
$\Lambda_{n} / (h)$ and 
$\Lambda_{p_1^{e_1}}/(h_1)  \times \ldots \times \Lambda_{p_k^{e_k}}/(h_k) $
are isomorphic, and since the quandle operations 
are defined using ring operations, they are isomorphic as quandles.
Then the result follows from Lemma~\ref{coeffcocylemma}.
\endproof

\begin{corollary}
Let $R_{n}$ be a dihedral quandle, where $n$ is a positive integer
with the prime decomposition $n=p_1^{e_1} \ldots p_k^{e_k}$,
for a positive integers $e_1, \ldots, e_k$ and $k$. 

Then  
the quandle 
$R_n$ is isomorphic to 
$R_{p_1^{e_1}} \times \ldots \times R_{p_k^{e_k}} $,
and each factor $R_{p_j^{e_j}}$ is inductively described
as 
an Alexander 
extension:
$R_{p_j^{d_j}} = AE( R_{p_j^{d_j -1}},\  R_{p}, \ \phi) . $
\end{corollary}

\begin{lemma}
Let $h \in \Lambda_p$ be a polynomial such that the coefficients of 
the highest and lowest degree terms are units in $\Z_p$. 
For any positive integer $m$, the Alexander quandle 
$E=\Lambda_p/(h^m)$ 
satisfies 
the conditions of Lemma~\ref{cocylemma},
with $X=\Lambda_p/(h^{m-1})$ and $A=\Lambda_p/(h)$.

Consequently, 
$$\Lambda_p/(h^m) =AE( \Lambda_p/(h^{m-1}),\ \Lambda_p/(h),\  \phi)$$
for some $\phi$ $\in$ $Z^2_{\rm TQ} ( \Lambda_p/(h^{m-1}); \Lambda_p/(h))$.
\end{lemma}
\proof
Assume that $h$ is a polynomial such that the lowest degree term is a 
non-zero constant,
and let $d=\mbox{deg}(h)$ be the degree of $h$. 

Define the map $f: E \rightarrow A \times X$ as follows.
Identify $\Lambda_p/(h^m)$ with $\Z_p [T]/ (h^m)$.
For a polynomial $L \in E$, 
write 
$$L=\sum_{j=0}^{m-1} A_j h^j= A_{m-1} h^{m-1} + \ldots + A_1 h + A_0 
= [A_{m-1}, A_{m-2}, \ldots, A_0], $$
 where $A_j \in \Z_p[T]$ has  degree less than $d$.
Let $$f(L)=( A_{m-1} \; \mbox{mod}\; (h),\ 
 \sum_{j=0}^{m-2}  A_j h^j\;  \mbox{mod}\; (h^{m-1})).$$
Denote $\overline{L}=  \sum_{j=0}^{m-2}  A_j h^j$,
 which is a well-defined polynomial, and denote 
$\tilde{L}=A_{m-1}  \mbox{mod}\; (h)$,
so that $f(L)=(\tilde{L}, \overline{L})$.

Let $s:  \Lambda_p/(h^{m-1}) \rightarrow  \Lambda_p/(h^{m}) $
be the set-theoretic section defined by 
$$s[ A_{m-2}, \ldots, A_0]= [0,  A_{m-2}, \ldots, A_0].$$ 
Let $s(\overline{L})*s(\overline{M})= [F_{m-1}, \ldots, F_0]$.

Let $L=\sum_j A_jh^j,  M=\sum_j B_jh^j  \in E$, then 
\begin{eqnarray*}
L*M &=& (T A_{m-1}+ (1-T)B_{m-1})h^{m-1} + s(\overline{L})*s(\overline{M}) \\
 &=& (\tilde{L} * \tilde{M} )h^{m-1} 
+\sum_{j=0}^{m-1} F_j  h^j, 
\end{eqnarray*}
and we have 
$$ f(L*M)=(\tilde{L} * \tilde{M} +F_{m-1} ,  \overline{L} *  \overline{M}). $$
Hence we have $\phi( \overline{L},  \overline{M})=F_{m-1}$. \endproof

\begin{theorem}\begin{sloppypar}
Let $\Lambda_p/ (h_1^{e_1} \ldots h_k^{e_k} )$ be 
an Alexander quandle, where $\{ h_1,$ $ \ldots,  h_k \} $ are polynomials 
such that the coefficients of the highest and lowest degree terms 
are units in $\Z_p$, and any pair of them is coprime,
where $k$ is a positive integer.
Then $\Lambda_p/ (h_1^{e_1} \ldots h_k^{e_k} )$ is isomorphic as 
quandles to 
$$ \Lambda_p/ (h_1^{e_1}) \times \ldots  \times  \Lambda_p/ (h_k^{e_k}), $$
and  
each 
factor is inductively described as Alexander 
extensions: 
$$ \Lambda_p/ (h_j^{d_j}) = AE( \Lambda_p/ (h_j^{d_j-1}), \ 
 \Lambda_p/ (h_j), \ \phi_j )$$
for some $\phi_j \in Z^2_{\rm TQ} ( \Lambda_p/ (h_j^{d_j-1});  \Lambda_p/ (h_j))$.\end{sloppypar}
\end{theorem}
\proof 
If $f, g \in \Lambda_p$ are coprime, then as 
$\Lambda$-modules, 
$ \Lambda_p/ (fg) $ is isomorphic to $ \Lambda_p/ (f) \times  \Lambda_p/ (g)$,
and the quandle structures on these 
$\Lambda$-modules 
 are defined by using the 
$\Lambda$-module 
structure 
so that they are isomorphic as quandles as well. 
The result, then, follows from the preceding lemma.
\endproof

\begin{example} \label{dihedex2} {\rm
For the extension $\Lambda_3/(T+1)^2 \!=\! AE(X, A, \phi')$ 
for $X\!=\!\Lambda_3/(T+1)=R_3= A$, 
computations that are similar to those in Example~\ref{dihedex}
gives the following $2$-cocycle $\phi'$:
$$ \phi'= 2 \chi_{0,1} + \chi_{0,2} + \chi_{1,0} + 2 \chi_{1,2} 
 + 2 \chi_{2,0} + \chi_{2,1} . $$} \end{example}

\begin{proposition} 
 $ \mbox{\rm Rank}\ H^2_{\rm TQ} (R_n;R_n) \geq 2$ if $n$ is odd. 
\end{proposition}
\proof 
Let $\phi, \phi'$ be cocycles defined by Alexander extensions
$\Lambda_{n^2} / (1+T)=R_{n^2}=AE(R_n, R_n, \phi)$
and $\Lambda_{n} / (1+T)^2=AE(R_n, R_n, \phi')$, respectively.
Let $x=(1,0)-(-1,0)$ and $y=(0,1)-(2,1)$, respectively. 
Then $x$ and $y$ are cycles, $x, y \in Z_2^{\rm TQ}(R_n;R_n)$, 
and satisfy 
$\phi (x)=-1$, $\phi(y)=0$, $\phi '(x)=2$, and $\phi' (y)=-2$.
\endproof

\begin{remark} {\rm
We conjecture that
$H^2_{\rm TQ}(A; A)$ has rank at least two, for any 
Alexander quandle of the form $A=\Lambda_n / (h)$,
where $n$ is a positive integer and $h$ is a polynomial with 
the leading and constant coefficient invertible.

} \end{remark}

\section{Cohomology with $H^1$ coefficients} \label{H1coeffsec}

\begin{sloppypar}
In this section we construct cocycles using one dimensional 
lower cocycles with $H^1$ coefficients.
Let $X$ be a finite quandle, and $A$ be a finite Alexander quandle.
Consider $\xi \in C^n_{\rm TQ}(X; H^1_{\rm TQ}(X;A)) $.
 For any $n$-tuple  $(x_1, \ldots, x_n)$
 of elements of $X$, $\xi (x_1, \ldots, x_n) \in  H^1_{\rm TQ}(X;A)= Z^1_{\rm TQ}(X;A)$.
Hence $\xi (x_1, \ldots, x_n)$ is a quandle homomorphism
$X \rightarrow A$, so that for any $x \in X$, we obtain 
 $\xi (x_1, \ldots, x_n) (x) \in A$. 
\end{sloppypar}

\begin{proposition} \label{h1cocyprop}
\begin{sloppypar}
Let $X$ be a finite quandle, and $A$ be a finite Alexander quandle.
If $\xi \in Z^n_{\rm TQ}(X; H^1_{\rm TQ}(X;A)) $ satisfies
\end{sloppypar}
$$ T \xi (x_1, \ldots, x_n) (x_{n+1}) 
= \xi (x_1 * x, \ldots, x_n*x)(x_{n+1}*x) $$
for any $x, x_1, \ldots, x_{n+1} \in X$, then 
$\psi \in  Z^{n+1}_{\rm TR}(X; A)$ where
$\psi$ is defined by 
 $\psi (x_1, \ldots, x_{n+1})=\xi  (x_1, \ldots, x_n) (x_{n+1})$.
\end{proposition}
\proof
We compute
\begin{eqnarray*}
\lefteqn{ (\delta \psi) (x_1, \ldots, x_{n+1}, x_{n+2})  } \\
 & = & \psi (\partial (x_1, \ldots, x_{n+1}, x_{n+2}) ) \\
 & = & (\delta \xi )  (x_1, \ldots, x_{n+1} ) ( x_{n+2}) \\
& & 
+ (-1)^n T \xi (x_1, \ldots, x_n)(x_{n+1}) \\ &&
-  (-1)^n \xi (x_1*x_{n+2},  \ldots, x_n*x_{n+2})
( x_{n+1}*x_{n+2}) 
\end{eqnarray*}
and the result follows by setting $x_{n+2}=x$. 
\endproof

\begin{example} {\rm \begin{sloppypar}
Let $X=A=R_3=\{ 0,1,2\}=\{a,b,c\}$. 
Let $\mu \in  C^1_{\rm TQ}(X; H^1_{\rm TQ}(X;A)) $.
The condition in Proposition~\ref{h1cocyprop} 
is written as \end{sloppypar}
$$ - \mu( x_1)(x_2)= \mu (x_1*x )(x_2*x) $$
for any $x,x_1,x_2 \in X=R_3$. 
We seek a $2$-cocycle $\phi (x_1, x_2)=\mu (x_1)(x_2) $,
 $\phi \in  Z^{2}_{\rm TQ}(X; A)$.
For the quandle cocycle condition ($\phi(x,x)=0$ for any $x \in R_3$), 
we assume  $\mu(x)(x)=0$. If $\mu(0)(1)=0$, then $\mu(0) \in H^1(R_3 ; R_3)$
is the constant homomorphism $\mu(0)(x)=0$ for any $x \in R_3$,
and a trivial $2$-cocycle $\phi$ results.
Hence we may assume that $\mu (0)(1)=1$ or $2$. 
Consider the case  $\mu (0)(1)=1$.
By the above formula, we have 
$$ \begin{array}{rrrrrrr}
\mu (0)(2) &= &\mu(0*0)(1*0) &= & -\mu (0)(1)& =& -1  \\
\mu (1)(2) &= &\mu(0*2)(2*2) &= & -\mu (0)(2)& =& 1 \\
\mu (1)(0) &=& \mu(1*1)(2*1) &=& -\mu (1)(2)& =& -1 \\
\mu (2)(0) &=& \mu(1*0)(0*0) &=& -\mu (1)(0)& =& 1 \\
\mu (2)(1) &=& \mu(2*2)(0*2) &=& -\mu (2)(0)& =& -1
\end{array}
$$
and  we obtain 
$$\phi= (\chi_{0,1} + \chi_{1,2}+ \chi_{2,0}) 
 + 2 ( \chi_{0,2} + \chi_{2,1}+ \chi_{1,0}) ,$$
which is the negative of the cocycle in Example~\ref{dihedex2}.
In fact, the case  $\mu (0)(1)=2$ yields the same cocycle as 
 Example~\ref{dihedex2}.

If we did not have this example in hand, then we are not yet able to 
conclude that the above obtained $\phi(x,y)=\mu(x)(y)$ is 
a cocycle, since we have not checked that 
$\mu \in Z^1_{\rm TQ}(R_3;R_3)$. However, from the above computations,
it is easily seen that for any $x$, $\mu(-)(x)$ is a quandle 
isomorphism on $R_3$,
as any permutation of the three elements is a quandle isomorphism.
Here, the second factor of $\mu$ is fixed and 
$\mu$ is regarded as a function with respect to the first factor.
 This fact of
 $\mu(-)(x)$ being isomorphisms is equivalent to
 $\mu \in Z^1_{\rm TQ}(R_3;R_3)$.

} \end{example}

\begin{example} \label{h3r3ex} {\rm \begin{sloppypar}
Again let $X=A=R_3$, and we construct a $3$-cocycle
$\theta \in Z^3_{\rm TQ}(R_3;R_3)$ by setting 
$\theta(x_1, x_2, x_3)=\phi(x_1, x_2)(x_3)$, where 
 $\phi$ $\in$ $C^2_{\rm TQ}$ $(X; H^1_{\rm TQ}(X;A))$. 
The condition in Proposition~\ref{h1cocyprop}
is written in this case as 
$ - \phi(x_1, x_2)(x_3)= \phi(x_1*x, x_2*x)(x_3*x) $ 
for any $x, x_1, x_2, x_3 \in R_3$. 
If $\phi(0,1)(0)=0$, then from  the quandle condition $\phi(0,1)(1)=0$,
we have the trivial homomorphism as $\phi (0,1)$, so that 
we assume  $\phi(0,1)(0)=1$ (the case  $\phi(0,1)(0)=-1=2$ yields the negative
of this case). 
For $\phi(0,1)$ to be an isomorphism of $R_3$, we have   $\phi(0,1)(2)=-1$.
Computations similar to the preceding example yield a $2$-cochain.
The computations are done by noticing the following sequence consisting
of actions by quandle elements from the right:
$$
\begin{array}{lllllllllll}
(0,1,0) & \stackrel{*0}{\rightarrow} & (0,2,0) 
& \stackrel{*1}{\rightarrow} & (2,0,2) 
&  \stackrel{*2}{\rightarrow} & (2,1,2) 
& \stackrel{*0}{\rightarrow} & (1,2,1) 
& \stackrel{*1}{\rightarrow} & (1,0,1) \\
(0,1,2)  & \stackrel{*0}{\rightarrow} & (0,2,1)
 & \stackrel{*1}{\rightarrow} & (2,0,1) 
 & \stackrel{*2}{\rightarrow} & (2,1,0) 
 & \stackrel{*0}{\rightarrow} & (1,2,0)
 & \stackrel{*1}{\rightarrow} & (1,0,2).
\end{array}
$$
This yields the cochain
\begin{eqnarray*}
\theta &=& (\chi_{0,1,0} + \chi_{2,0,2}+ \chi_{1,2,1} 
+ \chi_{0,2,1} + \chi_{2,1,0} + \chi_{1,0,2} ) \\
& & -( \chi_{0,2,0} + \chi_{2,1,2} +\chi_{1,0,1}
+ \chi_{0,1,2}+ \chi_{2,0,1}+\chi_{1,2,0} ). 
\end{eqnarray*}
It is checked that each $\phi(x,y)$ is in $H^1_{\rm TQ}(X;A)$,
 being a permutation. 
Now we check that $\phi(x,y) \in Z^2_{\rm TQ}(X; H^1_{\rm TQ}(X;A))$.
It is sufficient to prove that $\phi(x,y)(z)$ satisfies
 the $2$-cocycle condition for any $z \in R_3$.
{}From $\theta$ we have 
\begin{eqnarray*}
\phi(-,-)(0) &=& \chi_{0,1} +  \chi_{2,1} -  \chi_{0,2}- \chi_{1,2} \\
\phi(-,-)(1) &=& \chi_{1,2} +  \chi_{0,2} - \chi_{1,0} -  \chi_{2,0} \\
\phi(-,-)(2) &=& \chi_{2,0} + \chi_{1,0} -  \chi_{2,1} -  \chi_{0,1}.
\end{eqnarray*}
Let $\phi=\chi_{0,2}+\chi_{1,2} -\chi_{1,0} -\chi_{2,0} $
be the cocycle found in Example~\ref{R3ex}.
Note that $\delta \chi_0 = - \sum_{i \neq j} \chi_{i,j}$ where
the sum ranges over all pairs $(i,j) $, $i,j \in R_3$, such that $i\neq j$.
Then it is computed that 
$$ \phi(-,-)(0)=\phi - \delta \chi_0 , \quad
\phi(-,-)(1) = \phi , \quad 
\phi(-,-)(2)= \phi + \delta \chi_0,
$$
and we obtained  $\phi(x,y) \in Z^2_{\rm TQ}(X; H^1_{\rm TQ}(X;A))$.
Hence we constructed  
$\theta \in Z^3_{\rm TQ} (R_3;R_3)$ using Proposition~\ref{h1cocyprop},
from $Z^2_{\rm TQ} (R_3; H^1_{\rm TQ}(R_3;R_3)) $.
\end{sloppypar}
} \end{example}

\begin{proposition}
$H^3_{\rm TQ}(R_3;R_3) \neq 0$.
\end{proposition}
\proof Let $\theta \in Z^3_{\rm TQ} (R_3;R_3)$ 
be the cocycle obtained in Example~\ref{h3r3ex}. 
Let $c=(0,1,0)-(0,2,0) \in Z_2^{\rm TQ}(R_3;R_3)$.
It is easily computed that $c$ is indeed a  
$3$-cycle 
(see Example~\ref{R3ex}). Then it is evaluated that 
$\theta(c)=2 \neq 0$, hence $\theta \neq 0  \in H^3_{\rm TQ} (R_3;R_3)$.
\endproof

\section{Twisted cocycle knot invariants} \label{invsec}

We define the twisted cocycle knot invariant in this section.
First, we define the Alexander numbering for crossings.

\begin{figure}[ht!]
\begin{center}
\mbox{
\epsfxsize=2.5in
\epsfbox{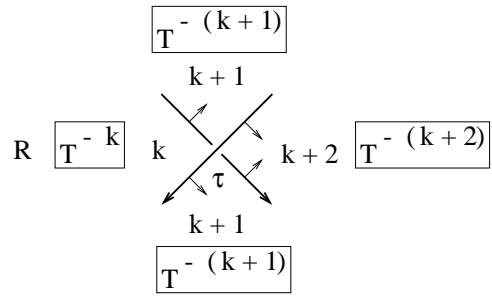} 
}
\end{center}
\caption{The Alexander numbering of a crossing }
\label{LXcrossing}
\end{figure}

Let $K$ be an oriented knot diagram with normals. 
Consider the underlying simple closed curve of $K$,
which is a generically immersed curve dividing the plane into regions,
and let $R$ be one of the regions.  
 Let $\alpha$ be an arc on the plane
from a point in the region at infinity to a point $R$ such that 
the  interior of $\alpha$ misses  all the 
crossing points of $K$ and intersects transversely in  finitely many 
points with the arcs of $K$. 
A classically known concept called {\it Alexander numbering} 
(see for example \cite{CS:book,CKamS:LX})
of $R$, denoted by $\LX(R)$, 
is defined as the number, counted with signs, 
of the  number of intersections between $\alpha$ and $K$.

More specifically, when $\alpha$ is traced from the region at infinity to 
$R$, and intersect at $p$ with $K$, if the normal to $K$ at $p$
is the same direction as $\alpha$, then $p$ contributes $+1$ to $\LX(R)$.
If the direction of $\alpha$ is the opposite to the normal, then 
its contribution is $-1$. The sum over all intersections
does not depend on the choice of $\alpha$.

In general, an Alexander numbering exists for an immersed 
curve in an orientable surface if and only if the curve
 represents a trivial $1$-dimensional class in the homology of the surface.

\begin{definition} {\rm
Let $K$ be an oriented knot diagram with normals. 
Let $\tau$ be a crossing.
There are four regions near $\tau$, 
and  the unique region
from which normals of over- and under-arcs point  
 is called the {\it source region} of $\tau$.

The {\it Alexander numbering} $\LX(\tau)$ of a crossing $\tau$ 
is defined to be 
 $\LX(R) $ 
 where $R$ is the source region of $\tau$. 
Compare with \cite{CKamS:LX}.

In other words, $\LX(\tau)$ is the number of intersections, counted with signs,
between an arc $\alpha$ from the region at infinity to
$\tau$ approaching from  the source region
of $\tau$.
 In  Fig.~\ref{LXcrossing}, the source region $R$ is the left-most region,
and 
the Alexander numbering of $R$ is $k$,
and so is the Alexander numbering of the crossing $\tau$.

} \end{definition}

Let a classical knot (or link) diagram $K$,  a finite quandle $X$,
a finite Alexander quandle $A$  
be given.  
A coloring 
of $K$ by $X$  also is given and 
is denoted by ${\cal C}$.

A {\it twisted (Boltzmann) weight}, $B_T(\tau, {\cal C})$,
at a  crossing $\tau$ is defined as follows.
Let  ${\cal  C}$ 
denote a coloring.
Let $r$ be the over-arc at $\tau$, and $r_1$, $r_2$ be 
under-arcs such that the 
normal
to $r$ points from $r_1$ to $r_2$.
Let $x={\cal C}(r_1)$ and $y={\cal C}(r)$. 
 Pick a 
quandle 
2-cocycle 
$\phi \in  Z^2_{\rm TQ}(X; A)$.
Then define 
$B_T(\tau, {\cal C})= [\phi(x,y)^{\epsilon (\tau)} ]^{T^{-\LX(\tau)}}$,
where 
$\epsilon (\tau)= 1$ or $-1$, if  the sign of $\tau$ 
is positive or negative, respectively.
Here, we use the multiplicative notation of elements of $A$, so that 
$\phi(x,y)^{-1}$ denotes the inverse of $\phi(x,y)$. 
Recall that $A$ admits an action by $\Z=\{ T^n \}$, 
and for $a \in A$, the action of $T$ on $a$ is denoted by $a^T$. 
To specify  the action by ${T^{-\LX(\tau)}}$ in the figures,
 each region $R$ with 
Alexander numbering
$\LX(R)=k$ is labeled  by the power $T^{-k}$ framed with a square,
as depicted in Fig.~\ref{LXcrossing}.

The {\it state-sum}, or a {\it partition function},
is the expression 
$$
\Phi (K) = \sum_{{\cal C}}  \prod_{\tau}  B_T( \tau, {\cal C}).
$$
The product is taken over all crossings of the given diagram,
and the sum is taken over all possible colorings.
The value of the weight  $B_T( \tau, {\cal C})$
is in the coefficient group $A$ written multiplicatively. 
Hence the value of the state-sum is in  the group ring ${\Z}[A]$.

\begin{figure}[ht!]
\begin{center}
\mbox{
\epsfxsize=3in
\epsfbox{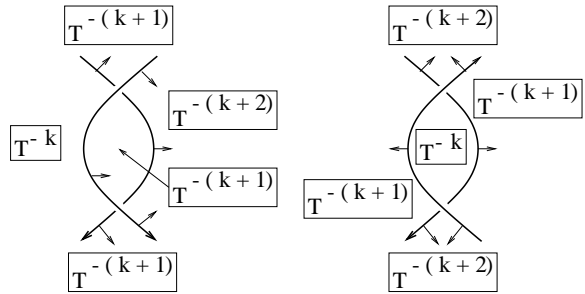} 
}
\end{center}
\caption{ Type II move and Alexander numbering }
\label{typeII} 
\end{figure}

\begin{figure}[ht!]
\begin{center}
\mbox{
\epsfxsize=3.8in
\epsfbox{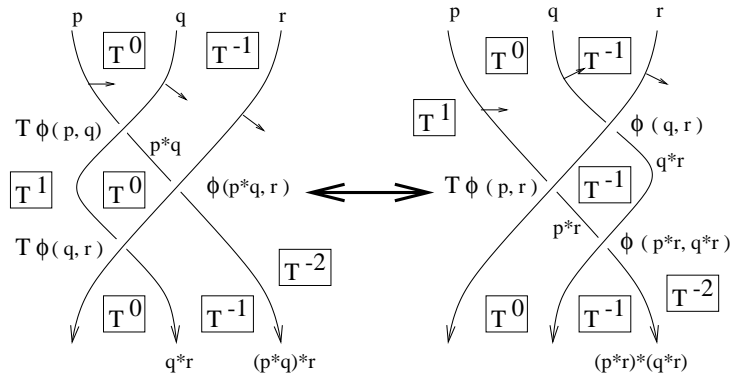} 
}
\end{center}
\caption{ Type III move and twisted $2$-cocycles }
\label{2tcocy} 
\end{figure}

\begin{theorem} 
The state-sum is well-defined.

More specifically, let  $ \Phi(K_1)$ and $ \Phi(K_2)$ 
be the state-sums  obtained from two diagrams of the same knot, then
 we have  $ \Phi(K_1) = \Phi(K_2)$. 
\end{theorem}
\proof
The invariance is proved by checking Reidemeister moves as follows.
Since the $2$-cocycle used satisfies $\phi(x, x)=1$ for any
$x \in X$, and the action of $T$ on the identity results in identity,
the type I Reidemeister move does not alter the state-sum.

For the type II move, we note that 
the  crossings involved in a type II move 
have opposite signs, and 
have the same Alexander numbering,
see Fig.~\ref{typeII} for typical situations 
(other cases can be checked similarly).
In both cases in the figure, all the crossings have
the same Alexander numbering $\LX(\tau)=k$, as seen from
the Alexander numberings of the adjacent regions specified in the figure
by square-framed labels. 
Hence the contribution to the state-sum of the pair of
crossings is of the form 
$ [  \phi(x,y)^{\epsilon} ]^{T^n} [\phi(x,y)^{-\epsilon}]^{T^n} $,
which is trivial.
Hence the state-sum is invariant under type II move.

Figure~\ref{2tcocy} depicts the situation for a type III move,
for  specific choices 
 of crossing information and orientations.
In this case, the left most crossings have the Alexander 
numbering $-1$
so that there is a $T$-factor in the Boltzmann weight,
and the right crossings, consequently, have numbering $0$, and do not 
have the $T$-factor. 
{}From the figure it is seen that the contributions to the state-sum,
in this case, is exactly the $2$-cocycle condition for the left and right hand 
side of the figure, and hence the state-sum remains unchanged.
In the figure, the $T$-action on cocycles is denoted 
in additive notation $T \phi(x,y)$ instead of
multiplicative notation $\phi(x,y)^T$, to match the $2$-cocycle condition
formulated in additive notation.
The other cases follow from combinations with type II moves,
see \cite{K&P,Turaev} and \cite{CJKLS} for more details.
\endproof

\begin{figure}[ht!]
\begin{center}
\mbox{
\epsfxsize=2in
\epsfbox{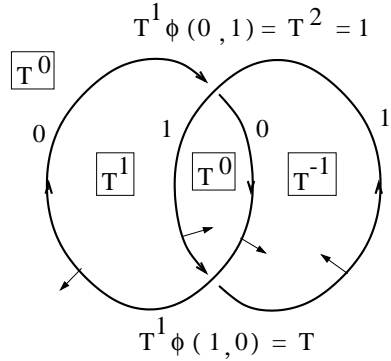} 
}
\end{center}
\caption{ Hopf link }
\label{hopf} 
\end{figure}

\begin{example} \label{hopfex} {\rm \begin{sloppypar}
Let $X=T_2$ (the trivial two element quandle) and $A=$ $\Z [T, T^{-1}]/(T^2-1)$.
$\phi=T \chi_{0,1} + \chi_{1,0}$
is a cocycle in $Z^2_{\rm TQ}(X;A)$. 
As an abelian group, $A$ is generated by $1$ and $T$, each 
denoted multiplicatively by $s$ and $t$, respectively. Thus 
any element of $A$ is written as $s^m t^n$ for integers $m,n$,
and the value of the invariant lies in 
$\Z [A]=\{ a+ b s^m t^n | a, b, m, n \in \Z \}$.
\end{sloppypar}

A coloring of a Hopf link $L$ and computations of weights are
depicted in Fig.~\ref{hopf}. This specific contribution to the state-sum
is $T+1$, or $st$. Note that both crossings
have the Alexander numbering $-1$, 
so the weight is multiplied by $T$.
By considering all possible colorings, 
we obtain $\Phi(L)= 2 + 2 st$. 

} \end{example}

For knots and links on compact surfaces defined up to 
Reidemeister moves, a similar invariants can be defined. 
There are two modifications that have to be made.

\noindent
(1) The regions divided by a given diagram have consistent 
colorings by powers of $T$.

\noindent
(2) Since there is no region at infinity, the choice of 
the ``base'' region must be considered.

Let $K$ be an oriented
 knot or link diagram on a compact oriented 
surface $F$.
Let $X$ be a finite quandle and $A$ be the coefficient group,
which is a $\Lambda=\Z [T, T^{-1}]$-module. 
Assume that $T^n=1$ for the action of $T$ on $A$. 
Let $\phi \in Z^2_{\rm TQ}(X;A)$. 

Let $R_i$, $i=0, 1, \ldots, n$, be the regions 
divided 
by
$K$, and call $R_0$ the {\it base} region. 
Define the {\it mod $p$ Alexander numbering} as before, 
except taking the values to be in $\Z_p$, where $p$ is a positive
integer.

If such a coloring of regions by $\Z_p$ is not possible, 
define $\Phi(K)=0$. 
Otherwise,  
we proceed as follows.
A  coloring 
${\cal C}$ of a knot diagram is defined similarly as before. 

\begin{sloppypar}
A {\it twisted (Boltzmann) weight}, $B_T(\tau, {\cal C})$,
at a  crossing $\tau$ is defined similarly by
$B_T(\tau, {\cal C})= [\phi(x,y)^{\epsilon (\tau)} ]^{T^{-\LX(\tau)}}$.
The {\it state-sum}, or a {\it partition function},  
is defined similarly by  
$
\Phi (K) = \sum_{{\cal C}}  \prod_{\tau}  B_T( \tau, {\cal C}).
$
\end{sloppypar}

To state the theorem, we need the following convention.
A typical element of  ${\Z}[A]$ is of the form $\sum_{i=1}^n x_i a_i$
for a positive integer $n$, where $x_i \in \Z$ and $a_i \in A$.
We define the {\it action} of $\Z=\langle T \rangle$ on  ${\Z}[A]$ by 
$ (\sum_{i=1}^n x_i a_i)^{T} = \sum_{i=1}^n x_i ( a_i)^{T}$. 
When a base region is replaced by another region, the state-sum 
changes by an action of $T^k$ for some integer $k$. 
Thus a proof similar to the planar diagram case implies the following 
generalization.

\begin{theorem} 
The state-sum is well-defined  up to the action of $\Z=\langle T \rangle$
for knots and links on surfaces.

More specifically, let  $ \Phi(K_1)$ and $ \Phi(K_2)$ 
be the state-sums  obtained from two diagrams of the same knot, then
for some integer $k$, we have  $ \Phi(K_1) = \Phi(K_2)^{T^k}$. 
\end{theorem}

\begin{remark} {\rm
For planar link diagrams, one could ``throw a string over the point 
at infinity,'' to shift  the Alexander numberings  by $\pm 1$.
The same change can be realized by Reidemeister moves.
This implies that the values of the invariant for planar link diagrams
are polynomials invariant under $T$-action. 
} \end{remark}

\begin{figure}[ht!]
\begin{center}
\mbox{
\epsfxsize=3in
\epsfbox{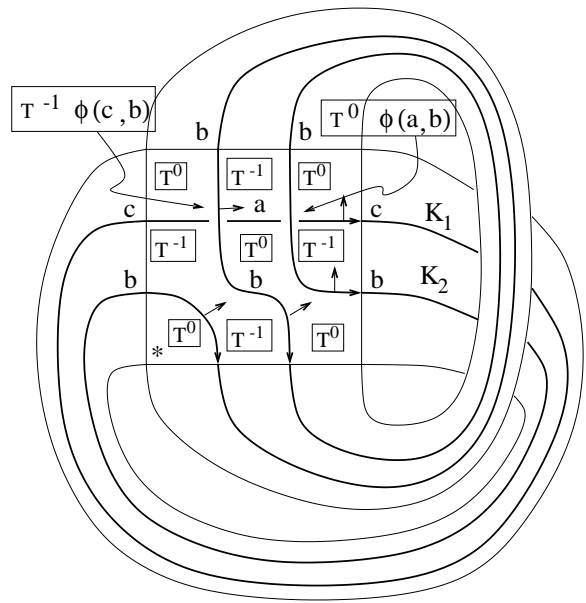} 
}
\end{center}
\caption{A link on a torus }
\label{t2cycle} 
\end{figure}

\begin{example} \label{torusex} {\rm 
A link $L$ on a torus is depicted in Fig.~\ref{t2cycle}. 
A coloring by $X=R_3=\{ a,b,c \}$ is given.
Note that the action of $T$ on $A=R_3$ satisfies
$T^2=1$, so that 
a  mod $2$ Alexander numbering is defined with $A=R_3$.
The base region is marked by $*$.
The base region has the Alexander numbering $0$, and is labeled with
 the $T$-term
$T^0$. The powers of $T$ that the other regions receive
are depicted in the figure.
Note that $T^{k+2}=T^k$, so that regions are labeled by
 either $T^0$ or $T^{-1}=T$.
The left/right sides and top/bottom sides of the middle square have identical 
colorings and numberings, respectively.
Thus these sides can be identified, as depicted, by bands, 
to obtain a punctured torus, and further the boundary can be 
capped off by a disk to obtain a torus. 
The contributions $\pm T^j \phi(x,y)$ to the Boltzmann weight 
of each crossing is indicated. 
For this specific coloring, the contribution is
$\phi(a,b)-\phi(c,b)$. 

{}From Example~\ref{dihedex2} we  have a $2$-cocycle 
$$ \phi'= 2 \chi_{0,1} + \chi_{0,2} + \chi_{1,0} + 2 \chi_{1,2} 
 + 2 \chi_{2,0} + \chi_{2,1} \in Z^2_{\rm TQ}(R_3;R_3). $$
With this cocycle, one computes that the 
invariant is $\Phi(L)=3+3t+3t^2$.
The action of $T$ on this element is 
$T \cdot (3+3t+3t^2)=3+3t^2+3t$ so that the action does not change this
element, 
and 
 the class of the polynomial $3+3t+3t^2$ under $T$-action
consists of a single element. 

It is seen that the invariant is trivial ($=9$)
if we use the cocycle $\phi$ in Example~\ref{dihedex}.

} \end{example}

\begin{figure}[ht!]
\begin{center}
\mbox{
\epsfxsize=3in
\epsfbox{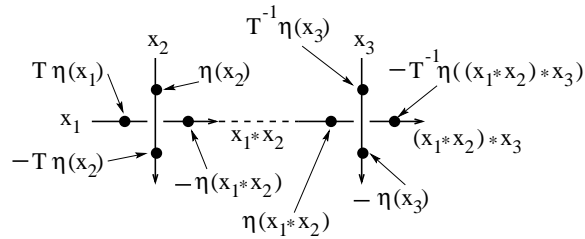} 
}
\end{center}
\caption{A coboundary defines the trivial invariant }
\label{cob} 
\end{figure}

\begin{proposition} \label{coblemma}
Let $X$ be a finite quandle, 
 and 
let 
$A$ be an Alexander quandle. 
Suppose $\phi \in Z^2_{\rm TQ}(X;A)$ is a coboundary: $\phi=\delta \eta$,
where $\eta \in  Z^1_{\rm TQ}(X;A)$.
Then the state-sum $\Phi(K)$ is a positive integer.
\end{proposition}
\proof
By assumption we have 
$$ \phi(x_1, x_2) = \delta \eta (x_1, x_2) = 
-T  \eta(x_2) +  T \eta(x_1) +  \eta(x_2) - \eta (x_1 * x_2). $$
For a given knot diagram $K$, remove a small neighborhood of each 
crossing, and let $\gamma_i$, $i=1, \ldots, m$, be the resulting arcs.
The end points of arcs are located near crossings, and depicted by
dots in Fig.~\ref{cob}. 
Assign each term of the above right-hand side to the end points as 
depicted in the left crossing of Fig.~\ref{cob}.
In the right of the figure, the 
situation at an adjacent crossing is also depicted. Note that 
the argument in $\eta$ 
coincides
with the color (a quandle element)
of the arc. Then it is seen that the terms assigned to 
the two end points of each arc are the same, with opposite signs
(as is seen from Fig.~\ref{cob}). 
Hence the contribution to the state-sum for any coloring is $1$,
and the state-sum is a positive integer (which is the number of 
colorings).
This argument is similar to the one given in \cite{CJKLS}.
\endproof

\begin{proposition} \label{obstvanishthm}
Let $\phi \in Z^2_{\rm TQ}(X;N)$ 
be 
an obstruction  $2$-cocycle,
where $X$ is a finite quandle and $N$ is an  Alexander quandle.
Then the state-sum invariant  $\Phi(K)$
defined from $\phi$ is a positive integer for any link diagram $K$
on the plane. 
\end{proposition}
\proof
We have an exact sequence 
$0  \rightarrow
 N \stackrel{i}{\rightarrow}  G \stackrel{p}{\rightarrow} A  \rightarrow 0$ 
of Alexander quandles, as in Theorem~\ref{obstthm},
and a section $s: A \rightarrow G$ with $ps=$id, $s(0)=0$.
By Relation~(\ref{2cocydefrel}), 
 for an obstruction cocycle 
$\phi$, we have 
$$ i\phi(x_1, x_2) =   Ts \eta(x_1) + (1-T) s\eta(x_2) - s\eta (x_1 * x_2). $$
Using $s \eta (x)$ instead of $\eta$ in the proof of the preceding 
Proposition, we obtain the result. 
Here, the fact that $K$ is a planar diagram is used in the 
step claiming that $ \pm T^k s \eta(x)$ assigned to endpoints of each arc 
cancel, since the $T$-factor $T^k$ matches on both endpoints of each arc.
More explanations on this point  are in order.
In the  
preceeding example 
of a link on a torus, the  Alexander 
numbering of regions satisfy $T^k=T^{k+2}$ since 
$T^2=1$ as 
an 
action on $N=R_3$, but the action of $T$ on the extension
$AE(R_3, R_3, \phi)$ does not
satisfy this relationship. 
Hence the terms
$T^k s (x)$ and $- T^{k-2} s (x)$ assigned to endpoints of a single 
arc do not cancel in the extension. In other words, 
in the preceding theorem, the cancelation was made in the coefficient ring,
but in this proof, the cancelations need to be done in the extension
via sections and inclusions, and the Alexander numbering of the regions
need to be consistent. The proof applies to such cases if 
the terms actually cancel, even if $K$ is non-planar.
\endproof

\begin{example} {\rm 
The link $L$ in Example~\ref{torusex}
has a non-trivial state-sum invariant with the cocycle 
 in Example~\ref{dihedex2}, which was obtained from 
a short exact sequence of Alexander quandles.
This is the case since, of course, $L$ is on a torus, and not on the plane.
} \end{example}

\begin{corollary}
Let $\phi \in Z^2_{\rm TQ}(X;A)$ 
be 
an obstruction  $2$-cocycle,
where $X$ and $A$ are finite Alexander quandles.
If the state-sum invariant   $\Phi(K)$
defined from $\phi$ is non-trivial (i.e., not a positive integer)
for a planar link diagram $K$, then  
the Alexander extension $AE(X, A, \phi)$ is not an Alexander quandle 
such that
$$ 0 {\rightarrow} A  \stackrel{i}{\rightarrow}A \times X = AE(X, A, \phi)
\stackrel{p}{\rightarrow}  X \rightarrow 0  $$ 
is a short exact sequence of $\Lambda$-modules where
 $i$ and $p$ are the natural maps
as in Remark~\ref{obstrem}. 
\end{corollary}
\proof
By Remark~\ref{obstrem}, if  $AE(X, A, \phi)$ is 
an Alexander quandle, then a short exact sequence
of Alexander quandles
$$ 0 \rightarrow A  \rightarrow A \times X = AE(X, A, \phi)
 \rightarrow  X \rightarrow 0  $$
defines an obstruction cocycle $\phi$. 
This contradicts the preceding Theorem.
\endproof

\begin{example} {\rm 
The $2$-cocycle $\phi$ $\in$ $Z^2_{\rm TQ}$  $(T_2;$ $\Z[T, T^{-1} ]/(T^2-1) )$ 
used in Example~\ref{hopfex} gave rise to a non-trivial value for 
a Hopf link. Hence 
$AE(T_2,$  $\Z[T, T^{-1} ]/(T^2-1) ,$ $ \phi)$ is not an Alexander quandle of 
the form 
stated 
in the preceeding Corollary.

For $A=T_2$, the cohomology theory is untwisted, 
and for $X=\Z_2[T, T^{-1}] / (T^2+T+1)$, it is known \cite{CJKS1} 
that $\phi = \sum_{a \neq b, a \neq T \neq b} \chi_{a,b}$ 
is a cocycle. With this cocycle, there are a number of classical
 knots in the table with non-trivial invariant. Hence
$AE(\Z_2[T, T^{-1}] / (T^2+T+1), T_2, \phi)$ is not an Alexander
quandle of the form 
stated 
in the preceeding Corollary.
} \end{example}

\begin{figure}[ht!]
\begin{center}
\mbox{
\epsfxsize=3in
\epsfbox{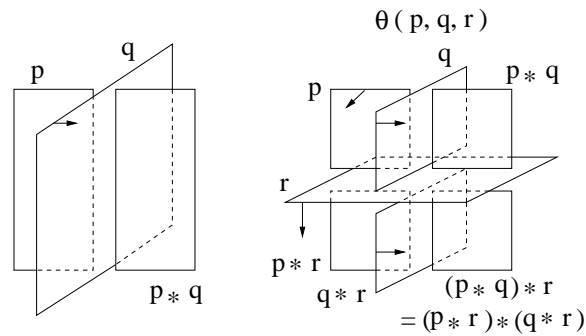} 
}
\end{center}
\caption{Colors at double curves and $3$-cocycle at a triple point }
\label{triplepoint} 
\end{figure}

The state-sum invariant is defined in an analogous way for 
oriented knotted surfaces
in $4$-space using their projections and diagrams in $3$-space. 
Specifically, the above steps can be repeated as follows,
for a fixed finite quandle $X$ and a knotted surface diagram $K$.

\begin{itemize}
\item
The diagrams consist of double curves and isolated branch and triple points
\cite{CS:book}. Along the double curves, 
the 
coloring rule is defined using normals
in the same way as classical case, as depicted in
 the left of Fig.~\ref{triplepoint}.

\item
The source region $R$ and the
Alexander numbering $\LX(\tau)=\LX(R)$ 
are  defined for a triple point
$\tau$  using normals.

\item 
A $3$-cocycle $\theta \in Z^3_{\rm TQ}(X;A)$, 
with the Alexander quandle coefficient $A$ is fixed, 
and assigned to a triple point as depicted in the right 
of Fig.~\ref{triplepoint}. In this figure, the triple point has the 
Alexander numbering $0$. 

\item
The sign $\epsilon(\tau)$ of a triple point $\tau$
is  defined \cite{CS:book}.

\item
For a coloring ${\cal C}$, the Boltzmann weight 
at a triple point $\tau$ is defined by 

\noindent
$B_T(\tau, {\cal C})= [\theta(x,y,z)^{\epsilon (\tau)} ]^{T^{-\LX(\tau)}}$. 

\item 
The state-sum is defined by 
$\sum_{{\cal C}}  \prod_{\tau}  B_T( \tau, {\cal C}).
$

 \end{itemize}

By checking the analogues of Reidemeister moves for knotted
surface diagrams, called Roseman moves, we obtain the following.

\begin{theorem}
The state-sum is well-defined for knotted surfaces,
and  is called the {\it twisted quandle cocycle invariant}
of knotted surfaces.
\end{theorem}

\begin{figure}[ht!]
\begin{center}
\mbox{
\epsfxsize=3in
\epsfbox{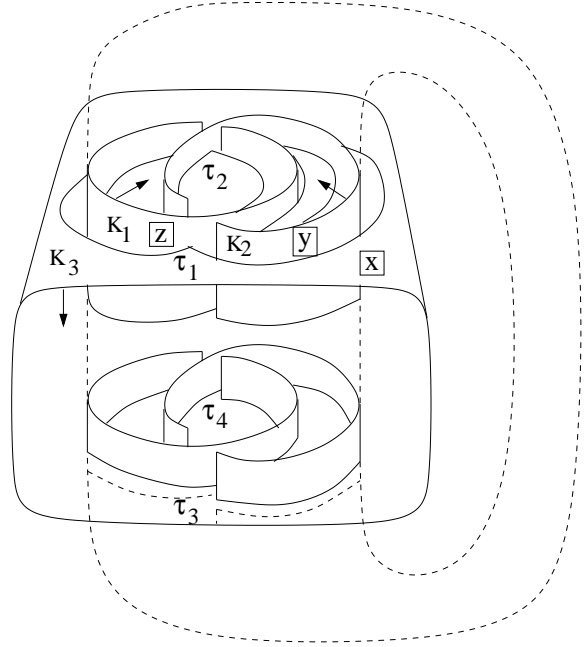} 
}
\end{center}
\caption{An analogue of Hopf link }
\label{bead} 
\end{figure}

\begin{example} \label{hopf4dex}   
{\rm 
Let $X=T_3=\{ 0,1,2\}$ (the trivial three element quandle) 
and $A=\Z [T, T^{-1}]/(T^2-1)$.
Recall that 
$\partial =(T-1) \partial_0$ as seen in Example~\ref{trivialex},
and $(T+1)(T-1)=0$ in $A$. 
It follows that 
$\theta=(T+1) \chi_{0,1,2}$
is a cocycle in $Z^3_{\rm TQ}(X;A)$ 
(in fact, this construction works in  Example~\ref{hopfex} as well). 
Denote the multiplicative generators of $A$ by $s$ and $t$, 
for additive generators $1$ and $T$, respectively.

In Fig.~\ref{bead}, an analogue of a Hopf link for surfaces in $4$-space,
$L=K_1 \cup K_2 \cup K_3$, is depicted.
Each component is standardly embedded in $4$-space, $K_1 \cup K_2$ is the spun
Hopf link with each component torus, and $K_3$ is a sphere
(in the figure, a large ``window'' is cut out from $K_3$ to show 
an inside view). 
The top horizontal sheet of $K_3$ is the bottom sheet
for the triple points $\tau_1$ and $\tau_2$
(that are positive triple points), and the bottom 
horizontal sheet of $K_3$ is the top sheet for $\tau_3$ and $\tau_4$
(that are negative triple points).
The orientation normals all point inside, so that all the triple points 
are negative, using the right-hand convention of the orientation
of the $3$-space. 
The source region is the region at infinity for all triple points,
so that the $T$-factor coming from the Alexander numbering is 
$T^0=1$ 
for all the triple points. 

The colors of relevant sheets are denoted by $x$, $y$, $z$, 
for sheets in $K_3$, $K_2$, and $K_1$, respectively,
as depicted. When trivial quandles are used, the colors 
depend only on the components. Hence the state-sum term 
is written by 
$$  \theta (x,y,z) \theta (x,z,y) \theta (y,z,x)^{-1}  \theta (z,y,x)^{-1} $$
where each term of $\theta$ coming from triple points 
$\tau_i$, $i=1,2,3,4$, respectively.

If the colors are given by $(x,y,z)=(0,1,2)$, 
$\theta(x,y,z)=T+1$ additively and $st$ multiplicatively, for example, and 
the above state-sum  term is equal to 
$ st $, since all the other $\theta$ terms are trivial.
The coloring $(x,y,z)=(0,2,1)$ also contributes $st$. 
The colorings $(x,y,z)=(2,0,1), (2,1,0)$
contributes $(st)^{-1}$.
All the other colorings contribute $1$, and the invariant 
is $\Phi(L)= 23+ 2 st + 2 (st)^{-1} $.

} \end{example}

\begin{figure}[ht!]
\begin{center}
\mbox{
\epsfxsize=3in
\epsfbox{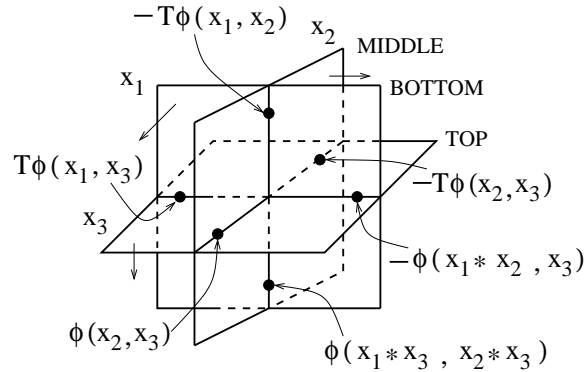} 
}
\end{center}
\caption{A coboundary at a triple point }
\label{cob3} 
\end{figure}

In fact, as in the classical case, the state-sum invariant is defined
modulo the action by $T$ for 
knotted surface diagrams in compact orientable $3$-manifolds,
up to Roseman moves. Such diagrams up to Roseman moves can be 
regarded as ambient isotopy  classes of embeddings of surfaces in the product 
space $M \times [0,1]$, where $M$ is a compact orientable $3$-manifold.

A similar argument to the proof of Proposition~\ref{coblemma} 
gives the following analogue, see Fig.~\ref{cob3}.
In this figure, a negative triple point is depicted, so that 
the terms are the negative of those that appear in $\delta \theta$. 
There is a diagram without branch point for orientable knotted surfaces
(see for example \cite{CS:cancel}), so that the terms assigned to 
the end points of double arcs cancel as in classical case, 
and we obtain the following.

\begin{proposition} \label{coblemma4D}
Let $X$ be a finite quandle, 
and 
let 
$A$ be an Alexander quandle. 
Suppose $\theta \in Z^3_{\rm TQ}(X;A)$ is a coboundary: $\theta=\delta \phi$,
where $\phi \in  Z^2_{\rm TQ}(X;A)$.
Then the state-sum $\Phi(K)$ for a knotted surface is a positive integer.
\end{proposition}

A similar argument to the proof of Proposition~\ref{coblemma4D}
and that of Theorem~\ref{obstvanishthm} 
can be applied to obtain the following.

\begin{proposition} \label{obst3cor}
Let $\theta \in Z^3_{\rm TQ}(X;N)$ be an  obstruction $3$-cocycle,
where $X$ is a finite   quandle and $A$ is an Alexander quandle. 
Then the state-sum invariant  $\Phi(K)$
defined from $\theta$ is a positive integer for any knotted surface 
diagram $K$ in Euclidean $3$-space ${\R}^3$. 
\end{proposition}

\begin{corollary}
Let $\theta \in Z^3_{\rm TQ}(X;N)$ an obstruction  $3$-cocycle,
where $X$ and $A$ are finite Alexander quandles.
If the state-sum invariant   $\Phi(K)$
defined from $\theta$ is non-trivial (i.e., not a positive integer)
for a knotted surface diagram $K$ in ${\R}^3$,
then $\theta$ is not an obstruction cocycle.
\end{corollary}

\begin{example} {\rm
By the preceding Corollary and Example~\ref{hopf4dex}, 
we find that
the cocycle in Example~\ref{hopf4dex}  
 is not an obstruction cocycle.
} \end{example}

\begin{figure}[ht!]
\begin{center}
\mbox{
\epsfxsize=4in
\epsfbox{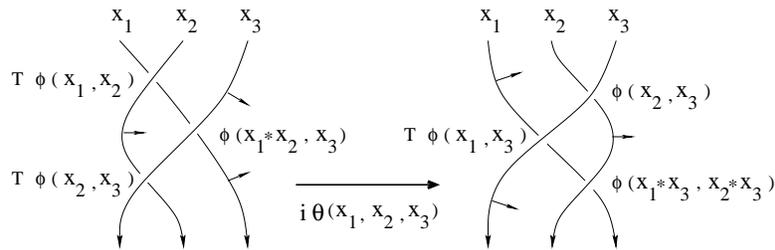} 
}
\end{center}
\caption{A $3$-cocycle assigned to a type III move }
\label{typeIII3cocy} 
\end{figure}

\begin{figure}[ht!]
\begin{center}
\mbox{
\epsfxsize=2.8in
\epsfbox{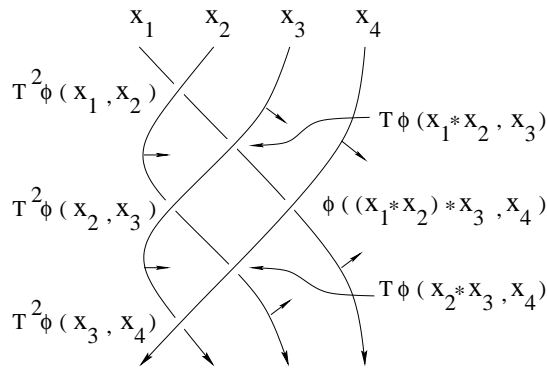} 
}
\end{center}
\caption{ The left-hand side of the $3$-cocycle condition}
\label{3cocyL} 
\end{figure}

\begin{figure}[ht!]
\begin{center}
\mbox{
\epsfxsize=4in
\epsfbox{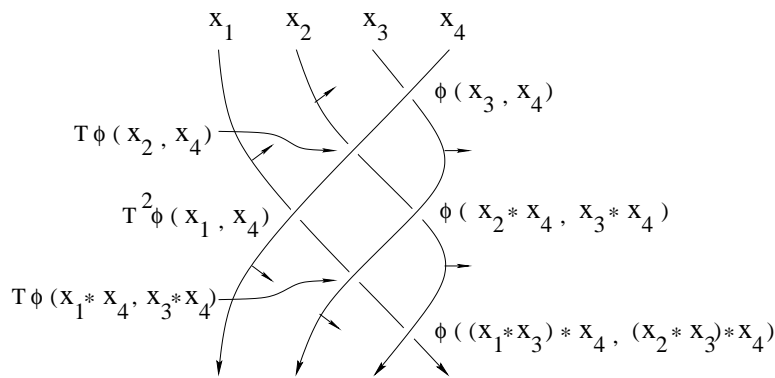} 
}
\end{center}
\caption{ The right-hand side of the $3$-cocycle condition}
\label{3cocyR} 
\end{figure}

\begin{remark} {\rm
As another application of colored knot diagrams, we exhibit a 
diagrammatic construction of the proof of Lemma~\ref{movielemma}.
Diagrammatic methods in cohomology theory, such as Hochschild cohomology,
are found, for example, in \cite{MS}.

In the state-sum invariant, 
 a $3$-cocycle is assigned to  a triple point as a Boltzmann weight.
When a height function in $3$-space is chosen, a triple point is described by 
the Reidemeister type III move. Cross sections of three sheets 
at a triple point by planes normal to 
the chosen height function give rise to a move among three strings, 
and the move is exactly the type III move. See \cite{CS:book} for more details.
In Fig.~\ref{typeIII3cocy}, the type III move as
such a movie description of a colored triple
point is depicted. 
In this movie, we color the diagrams by quandle elements, assign $2$-cocycles
to crossings, assign $3$-cocycles to type III move performed,
and the convention of these assignments is depicted in 
 Fig.~\ref{typeIII3cocy}.

In Figs.~\ref{3cocyL} and \ref{3cocyR},
 diagrams involving four strings are  depicted. 
These are cross sections of three coordinate planes in $3$-space
plus another plane in general position with the coordinate planes.
 See \cite{CS:book} for more details.
The colorings by quandle elements and  $2$-cocycles 
 are also depicted. 
Note that the $2$-cocycles depicted in 
Fig.~\ref{3cocyL} are exactly  the first expression  of 
in the proof of Lemma~\ref{movielemma}, 
and  those in  Fig.~\ref{3cocyR} are   the last expression,
respectively.

There are two distinct  sequences of type III moves 
that change Fig.~\ref{3cocyL} to Fig.~\ref{3cocyR}.
Each type III move gives rise to a $3$-cocycle via 
the convention established in Fig.~\ref{typeIII3cocy}.
It is seen that the two sequences of $3$-cocycles 
corresponding to two sequences of type III moves are identical to 
the sequences of equalities in the proof of Lemma~\ref{movielemma}. 
Once the direct correspondence is made, the computations follows 
from these diagrams 
automatically. 
} \end{remark}

\Addresses\recd

\end{document}